\documentclass[12pt,psamsfonts]{amsart}

\usepackage{amsmath,amssymb}
\usepackage{graphicx}
\usepackage{color}
\usepackage{psfrag} 
\usepackage{overpic} 
\usepackage{array} 
\usepackage{subfigure}
 
\newtheorem{theorem}{Theorem}
\newtheorem{proposition}[theorem]{Proposition}

\newtheorem{lemma}[theorem]{Lemma}
\newtheorem {remark}[theorem]{Remark}
\newtheorem{definition}[theorem]{Definition}

\title[Nilpotent Centers in $\mathbb{R}^3$]{Nilpotent Centers in $\mathbb{R}^3$}
\author[L. Queiroz and C. Pessoa]{}


\begin{document}
 \maketitle

\centerline{\scshape  Claudio Pessoa,  \; Lucas Queiroz}
\medskip

{\footnotesize \centerline{Universidade Estadual Paulista (UNESP), Instituto de Bioci\^encias Letras e Ci\^encias Exatas,} \centerline{R. Cristov\~ao Colombo, 2265, 15.054-000, S. J. Rio Preto, SP, Brasil }
\centerline{\email{c.pessoa@unesp.br} and \email{lucas.queiroz@unesp.br}}}

\medskip

\bigskip

\begin{quote}{\normalfont\fontsize{8}{10}\selectfont
{\bfseries Abstract.}
Consider analytical three-dimensional differential systems having a singular point at the origin such that its linear part is $y\partial_x-\lambda z\partial_z$ for some $\lambda\neq 0$. The restriction of such systems to a Center Manifold has a nilpotent singular point at the origin. We study the formal integrability and the center problem for those types of singular points in the monodromic case. Our approach do not require polynomial approximations of the Center Manifold in order to study the center problem. As a byproduct, we obtain some useful results for planar $C^r$ systems having a nilpotent singularity. We conclude the work solving the Nilpotent Center Problem for the Generalized Lorenz system and the Hide-Skeldon-Acheson dynamo system.
\par}
\end{quote}

\section{Introduction}

Consider a planar vector field $X$ defined in a open set $U\subset\mathbb{R}^2$. We say that a singular point $p$ is a nilpotent singular point of $X$ if the Jacobian matrix $DX(p)$ has two zero eigenvalues but is not the null matrix. In other words, $p$ is a nilpotent singular point of vector field $X$ if there exists a linear change of coordinates such that, in those new coordinates, $DX(p) =\begin{tiny}
	\left(\begin{array}{lr}
		0 & 1 \\
		0 & 0
	\end{array}\right)
\end{tiny}$. We may assume, by means of translation, that $p=0$. Thus, the associated differential system can be written as:
\begin{equation}\label{eqPlanar}
	\begin{array}{lr}
		\dot{x}=y+P(x,y),\\
		\dot{y}=Q(x,y),
	\end{array}
\end{equation}
where $j^1P(0)=j^1Q(0)=0$ (here $j^lF(0)$ denotes de $l$-jet of a function $F$ at $0$, i.e. the Taylor polynomial of degree $l$ at $0$). System \eqref{eqPlanar} is widely studied in the literature. There are several possibilities of phase portraits of this system in a neighborhood of the origin. In fact, almost all possibilities are classified, except the monodromic case (see  \cite[Theorem 3.5, p. 116]{Llibre} and also \cite{Andreev,Andronov}). 

%

In the case which \eqref{eqPlanar} is analytical, the only monodromic possibilities for the origin to be are the nilpotent center or the nilpotent focus. The problem of distinguishing between these characterizes the so-called \emph{Nilpotent Center Problem} \cite{AlvarezGasull1,Zoladek}. For analytic system \eqref{eqPlanar}, the following result (see \cite[Theorem 3.5, p. 116]{Llibre} and \cite[Theorem 66, p. 357]{Andronov}) provides a monodromy criterion. 

\begin{theorem}[Andreev's Theorem]\label{TeoAndreev}
	Let $X$ be the vector field associated to analytic system \eqref{eqPlanar} and the origin be an isolated singular point. Let $y=F(x)$ be the solution of the equation $y+P(x,y)=0$ in a neighborhood of $(0,0)$ and consider $f(x)=Q(x,F(x))$ and $\Phi(x)=\mbox{\rm div}X\vert_{(x,F(x))}$. Then, we can write
	$$f(x)=ax^{\alpha}+O(x^{\alpha+1}),$$
	$$\Phi(x)=bx^{\beta}+O(x^{\beta+1}).$$
	The origin is monodromic if and only if $a<0,\alpha=2n-1$ and one of the following conditions holds:
	\begin{itemize}
		\item[i)]$\beta>n-1$ or $\Phi\equiv 0$;
		\item[ii)]$\beta=n-1$ and $b^2+4an<0$.
	\end{itemize}
\end{theorem}

In this work, we study analytic vector fields $X$ in $\mathbb{R}^3$ having a singular point $p$ such that $DX(p)$ has the following Jordan canonical form
$$\begin{small}
	\left(\begin{array}{lcr}
		0 & 1 & 0 \\
		0 & 0 & 0 \\
		0 & 0 & -\lambda
	\end{array}
	\right),
\end{small}$$
where $\lambda\neq 0$. We can assume that the singular point is located at the origin. The associated differential system can be written, by a linear change of variables, in the following form:
\begin{equation}\label{eq1}
	\begin{array}{lr}
		\dot{x}=y+P(x,y,z),\\
		\dot{y}=\hspace{0.7cm}Q(x,y,z),\\
		\dot{z}=-\lambda z + R(x,y,z),
	\end{array}
\end{equation}
where $P,Q,R$ are analytic functions and $j^1P(0)=j^1Q(0)=j^1R(0)=0$. For such systems, the following theorem holds.

\begin{theorem}[Center Manifold Theorem]\label{TeoCenterManifold} Consider system \eqref{eq1}. There exists an invariant bidimensional $C^r$-manifold tangent to the $xy$-plane at the origin for every $r\geqslant 1$.
\end{theorem}

This theorem is proven in \cite{Kelley} and the manifolds described in its statement are called \emph{Center Manifolds}. An extensive study on the center manifolds can be found in \cite{Sijbrand}. We remark that neither the analiticity nor the uniqueness of the center manifold is guaranteed by Theorem \ref{TeoCenterManifold}. 

Since every center manifold is tangent to the $xy$-plane at the origin, it has a local parametrization as $z=h(x,y)$, where $h:U\subset\mathbb{R}^2\to\mathbb{R}$ is a $C^r$-function with $j^1h(0,0)=0$. Substituting $z=h(x,y)$ in \eqref{eq1}, we obtain a restricted bidimensional differential system given by
\begin{equation}\label{eq1restricted}
	\begin{array}{lr}
		\dot{x}=&y+P(x,y,h(x,y)),\\
		\dot{y}=&Q(x,y,h(x,y)).
	\end{array}
\end{equation}
This restricted system has a nilpotent singular point at the origin. Henceforth we say that a three-dimensional analytical vector field has a \emph{nilpotent singular point} if its associated differential system can be written as \eqref{eq1}.

The goal of this work is to study the Nilpotent Center Problem for system \eqref{eq1}, that is, to determine whether the origin is a nilpotent center or a nilpotent focus on a center manifold. This problem is barely explored in the literature and, as far as we know, there are only two works on this subject \cite{SongWangFeng, Wang2}. In these papers, the authors study particular families of system \eqref{eq1} and the main technique used to obtain the results consists in restricting the system to a center manifold using a truncated expression of its parametrization. This method is not very practical for computational purposes. Our approach seeks to study the nilpotent center problem without going through the restriction of the system to a center manifold.

The structure of our paper is as follows. In Section \ref{SecR2} we exhibit the fundamental results for planar systems having nilpotent singularities, starting with the problem of determine when the singular point is monodromic. We present a version of Andreev's Theorem for the $C^r$ case which allow us to identify among systems \eqref{eq1} those for which the origin is monodromic on a center manifold. The number $n$ in Theorem \ref{TeoAndreev} is called \emph{Andreev number} and is invariant by orbital equivalence. We characterize the systems \eqref{eq1} with Andreev number $2$. For nilpotent singular points with odd Andreev number, we prove that the monodromy condition $\beta=n-1$ implies that the singularity is a focus for system \eqref{eqPlanar}. In this development we also prove that each one of the monodromy conditions (i) and (ii) from Theorem \ref{TeoAndreev} are invariant by local diffeomorphisms. 

In Section \ref{SecLinearOperators} we work with some auxiliary linear operators which frequently show up in our investigation. Using those linear operators, we obtain a formal normal form for system \eqref{eq1}, which is presented in Section \ref{SecNormalForms}.

In Section \ref{SecFI} we study the relationship between formal integrability and the Nilpotent Center Problem for three-dimensional system \eqref{eq1}. We prove that, for the vector field $X$ associated to system \eqref{eq1}, there exists a formal series $H(x,y,z)=y^2+\sum_{n\geqslant 3}H_n(x,y,z)$ such that ${XH=\sum_{n\geqslant 4}\omega_nx^n}$. Moreover, if the first non-zero $\omega_n$ has even index $n$, then the origin of \eqref{eq1} cannot be a center. This result provides a method for the study of the Nilpotent Center Problem which does not require a polynomial approximation of the Center Manifold. Using the invariance of the monodromy conditions, proven in Section \ref{SecR2}, we were able to prove that system \eqref{eqPlanar} satisfying the monodromy condition $\beta=n-1$ is not formally integrable. As far as we know, this result is not found in the literature. We generalize this result to $\mathbb{R}^3$. Moreover, for system \eqref{eq1} under monodromy condition $\beta=n-1$ with $n$ odd, we prove that the origin cannot be a center on a center manifold. Furthermore, we classify the normal forms for systems \eqref{eq1}, having the origin as a monodromic singular point, that are formally integrable.

Finally, in Section \ref{SecApplications}, we solve the Nilpotent Center Problem for some three-dimensional systems, including the Generalized Lorenz system and the Hide-Skeldon-Acheson dynamo system.

\section{Fundamental results in $\mathbb{R}^2$}\label{SecR2}

\subsection{Monodromy Criterion}

The lack of analyticity of the restricted system \eqref{eq1restricted} apparently poses an obstacle for us to use Theorem \ref{TeoAndreev} to study the monodromy of the singularity on a center manifold. However, the proof presented in \cite{AlvarezGasull1} uses the concept of the Poincaré index which is not dependent on the analyticity of the vector field. The other important step of the proof uses the \emph{generalized polar coordinates}, first introduced by Lyapunov \cite{LyapunovP} (see also \cite{GasullTorre}), which also do not depend on the analyticity. Thus, by supposing $j^rf(0)\neq 0$, the same arguments can be made and we have a less restrictive version of Theorem \ref{TeoAndreev}.

\begin{theorem}[$C^r$-Andreev's Theorem]\label{TeoAndreevRestrito}
	Let $X$ be the vector field associated to the $C^r$-system, $r\geqslant 3$, given by
	\begin{equation}\label{eqPlanarCr}
		\begin{array}{lr}
			\dot{x}=y+X_2(x,y),\\
			\dot{y}=Y_2(x,y),
		\end{array}
	\end{equation} 	
	where $X_2,Y_2\in C^r$, $j^1X_2(0)=j^1Y_2(0)=0$ and such that the origin is an isolated singular point. Let $y=F(x)$ be the solution of the equation $y+X_2(x,y)=0$ through $(0,0)$  and consider $f(x)=Y_2(x,F(x))$ and $\Phi(x)=\mbox{\rm div}X\vert_{(x,F(x))}$. We can write
	$$f(x)=ax^{\alpha}+O(x^{\alpha+1}),$$
	$$\Phi(x)=bx^{\beta}+O(x^{\beta+1}).$$
	for $\alpha<r$. Suppose that $a\neq 0$, then the origin is monodromic if and only if $a<0, \alpha=2n-1$ and one of the following conditions holds:
	\begin{itemize}
		\item[i)]$\beta>n-1$ or $j^r\Phi(0)\equiv 0$;
		\item[ii)]$\beta=n-1$ and $b^2+4an<0$;
	\end{itemize}
\end{theorem}

For the sake of completeness, we present the proof of this theorem in the appendix. The above theorem provides a tool to study the monodromy of the origin for system \eqref{eq1} on a center manifold. In order to perform this study, we apply Theorem \ref{TeoAndreevRestrito} to the restricted system \eqref{eq1restricted} and check if the function $f$ is not flat, since for any $r\geqslant 3$, there exists a center manifold such that \eqref{eq1restricted} is a $C^r$-system.

The positive integer $n$ in the statements of Theorems \ref{TeoAndreev} and \ref{TeoAndreevRestrito} plays an important role in the study of nilpotent monodromic singular points. So we define the \emph{Andreev number} of a nilpotent singular point by the number $n$ in function $f(x)=ax^{2n-1}+O(x^{2n})$ described in Theorems \ref{TeoAndreev} and \ref{TeoAndreevRestrito}.
In \cite{GarciaFiI} the author proves that the Andreev number is invariant by analytical and formal orbital equivalence, i.e. via analytical and formal diffeomorphisms and time rescalings.

The following result is essentially a shortcut to verify monodromy for a given bidimensional system. 


\begin{lemma}[Monodromy criterion]\label{LemaMonodromiaDelta}
	Consider a $C^r$-system \eqref{eqPlanar}, the functions $F(x),f(x),\Phi(x)$ given in Theorem \ref{TeoAndreevRestrito}. Let $\tilde{a}$ be the coefficient of $x^{2n-1}$ in the power series expansion of $f(x)$, $\tilde{b}$ the coefficient of $x^{n-1}$ in the power series expansion of $\Phi(x)$ and $\Delta=\tilde{b}^2+4\tilde{a}n$. Suppose $j^{2n-2}f(0)=0$ and $j^{n-2}\Phi(0)=0$, then the singular point is monodromic with Andreev number $n$ if and only if $\Delta<0$. 
\end{lemma}
\noindent\textbf{Proof: }Note that $\Delta<0$ gives us immediately the first monodromy condition, i.e. $\alpha=2n-1$ and $a=\tilde{a}<0$. If $\tilde{b}=0$, then either $j^r\Phi(0)\equiv 0$ or $\Phi(x)\in O(x^n)$ and the monodromy condition (i) in Theorem \ref{TeoAndreevRestrito} is satisfied. If $\tilde{b}\neq 0$, the condition (ii) in Theorem \ref{TeoAndreevRestrito} is satisfied. Thus, the singular point is monodromic with Andreev number $n$. Conversely, suppose the singular point is monodromic with Andreev number $n$. Then $\tilde{a}<0$ and one of the monodromy conditions (i) or (ii) in Theorem \ref{TeoAndreevRestrito} must hold. If (i) holds, $\tilde{b}$ is null and $\tilde{a}<0$ implies $\Delta<0$. On the other hand, if (ii) is satisfied then $\Delta$ must be negative. Therefore $\Delta<0$ regardless the monodromy condition.\qed

\medskip

Armed with the above results, we can procedurally apply an algebraic algorithm to obtain conditions on the parameters of system \eqref{eq1} for the origin to be monodromic on a center manifold. This algorithm is described in the next section.

\subsection{The algorithm to detect monodromy}

Consider the following representation of analytic system \eqref{eq1}:
\begin{equation}\label{eq2}
	\begin{array}{lr}
		\dot{x}=&y+\sum\limits_{j+k+l\geqslant 2}a_{jkl}x^jy^kz^l,\\
		\dot{y}=&\sum\limits_{j+k+l\geqslant 2}b_{jkl}x^jy^kz^l,\\
		\dot{z}=&-\lambda z+\sum\limits_{j+k+l\geqslant 
			2}c_{jkl}x^jy^kz^l.
	\end{array}
\end{equation}
By the Center Manifold Theorem \ref{TeoCenterManifold}, any center manifold for \eqref{eq2} has a local parametrization $z=h(x,y)$ for which $h\in C^r$ and $j^1h(0)=0$. Being an invariant manifold, it must satisfy the following equation:
\begin{tiny}
\begin{equation}\label{eqAlgoWc}
	\dfrac{\partial h}{\partial x}(x,y)(y+P(x,y,h(x,y)))+\dfrac{\partial h}{\partial y}(x,y)(Q(x,y,h(x,y)))+\lambda h(x,y)-R(x,y,h(x,y)) = 0.
\end{equation}
\end{tiny}
This equation allows us to obtain the Taylor polynomial of $h(x,y)$ to any finite order $m$. In order to achieve this, we write $h(x,y)$ as a formal series $h(x,y)=\sum_{j+k\geqslant 2}h_{jk}x^jy^k$ and substitute in \eqref{eqAlgoWc}. We expand the left-hand side of this equation in a power series of $x,y$ and equate the coefficients to zero to obtain the coefficients $h_{jk}$ given in terms of the coefficients $a_{jkl},b_{jkl},c_{jkl}$ from \eqref{eq2}.

More precisely, we start by equating the coefficients of the homogeneous terms of degree $2$, which is the lowest, to zero. Then, we determine $h_{20},h_{11},h_{02}$. The coefficients $h_{jk}$ of homogeneous degree $3$, i.e. such that $j+k=3$, are given by the coefficients $a_{jkl},b_{jkl},c_{jkl}$ of system \eqref{eq2} and the coefficients $h_{jk}$ of the previous homogeneous degree. We repeat this iterative process obtaining successively the coefficients $h_{jk}$ of a homogeneous part once having computed the coefficients $h_{jk}$ of the previous ones.

Once the coefficients $h_{jk}$ are obtained up to a large enough homogeneous degree $m$, we have successfully retrieved the Taylor polynomial $j^mh(0)=\sum_{j+k=2}^{m}h_{jk}x^jy^k$. It is well known that any $C^r$ center manifold, for $r\geqslant m$, has the same $m$-jet (see, for instance, \cite{Sijbrand}). Substituting $h(x,y)=j^mh(0)+O(\Vert x,y\Vert^{m+1})$ in \eqref{eq2} yields a restricted $C^r$-system \eqref{eq1restricted} for which the $m$-jet of each component of the associated vector field is known.

We can proceed analogously to retrieve the $m$-jet of function $F(x)$ from Theorem \ref{TeoAndreevRestrito}. It is sufficient to consider the equation
\begin{equation*}
	F(x)+P(x,F(x),h(x,F(x)))=0.
\end{equation*}
Consequently, determining $j^mF(0)$ also determines $j^mf(0)$ and $j^m\Phi(0)$. Hence, if there exists an $m\in\mathbb{N}$ such that $j^mf(0)\neq 0$, we can apply this process to study the monodromy using Theorem \ref{TeoAndreevRestrito}.

Applying the algorithm described above to $m=3$, we obtain:
$$j^2h(0)=\frac{c_{200}}{\lambda}x^2+\frac{(\lambda c_{110}-2c_{200})}{\lambda^2}xy+\frac{(2c_{200}-\lambda c_{110}+c_{020}\lambda^2)}{\lambda^3}y^2.$$
Which yields the following:
\begin{equation}\label{eqfPhi}
	\begin{array}{lcr}
		j^3f(0)=b_{{200}}{x}^{2}-{\dfrac { \left( b_{{110}}a_{{200}}\lambda-b_{{300}}\lambda-b_{{101}}c_{{200}} \right) {x}^{3}}{\lambda}},\\
		j^1\Phi(0)=(2a_{200}+b_{110})x.
	\end{array}
\end{equation}
The full expressions of $j^3h(0)$ and $j^3\Phi(0)$ are too extensive, so we omit them here. With the information given by \eqref{eqfPhi}, it is possible to conclude, by Theorem \ref{TeoAndreevRestrito}, that monodromy can only occur if ${b_{200}=0}$. Furthermore, we have the necessary tools to classify the first family of systems having a monodromic singular point at the origin on a center manifold.

\begin{proposition}\label{PropoMono2}
	The origin is a nilpotent monodromic singular point with Andreev number $2$ on a center manifold of system \eqref{eq2} if and only if $b_{200}=0$ and
	\begin{equation}\label{ineqMono2}
		\dfrac{b_{{101}}c_{{200}}}{\lambda}<-\frac{(2a_{200}-b_{110})^2}{8}-b_{300}.
	\end{equation}
	Moreover, if ${2a_{200}+b_{110}}\neq 0$, the restricted system satisfies the monodromy condition $\beta=n-1$ in Theorem \ref{TeoAndreevRestrito}. 
\end{proposition}

\noindent\textbf{Proof: }The origin is monodromic only if $b_{200}=0$. Thus, by Lemma \ref{LemaMonodromiaDelta}, the necessary and sufficient condition for the Andreev number to be $2$ is $\Delta<0$, which by \eqref{eqfPhi} becomes:
$$\frac { \left( 4\,{a_{{200}}}^{2}-4\,b_{{110}}a_{{200}}+{b_{{110}}}^{2}+8\,b_{{300}} \right) \lambda+8\,b_{{101}}c_{{200}}}{
	\lambda}<0,$$
or equivalently
$$\dfrac{b_{{101}}c_{{200}}}{\lambda}<-\frac{(2a_{200}-b_{110})^2}{8}-b_{300}.$$
The second statement also follows from \eqref{eqfPhi}. \qed
\medskip


%
%

\subsection{Some results for planar $C^r$ systems with nilpotent singular points}

We now prove some results for planar $C^r$-systems which are useful for the Nilpotent Center Problem in $\mathbb{R}^3$. First, we need the following definitions.

\begin{definition}
	A polynomial $p\in\mathbb{R}[x,y]$ is a \emph{$(t_1,t_2)$-quasi-\-homogeneous polynomial of weighted degree $k$} if $p(\lambda^{t_1} x, \lambda^{t_2} y)=\lambda^kp(x,y)$. A general expression for such polynomials is $p(x,y)=\sum_{t_1i+t_2j=k}a_{ij}x^iy^j$ where $a_{ij}\in\mathbb{R}$. The vector space of all $(t_1,t_2)$-quasi-\-homogeneous polynomial of weighted degree $k$ is denoted by $\mathcal{P}^{(t_1,t_2)}_{k}$.
\end{definition}

\begin{proposition}\label{PropoPreformMono}
	Consider a $C^r$-system \eqref{eqPlanar}. Let $n$ be a positive integer such that $2n-1<r$. The origin is a monodromic isolated singular point with Andreev number $n$ if and only if there exists a local analytical change of variables that transforms \eqref{eqPlanar} into the following form
	\begin{equation}\label{eq1nHomo}
		\begin{array}{lcr}
			\dot{x}=y+\mu x^n+\displaystyle\sum_{i+nj\geqslant n+1}\tilde{a}_{ij}x^iy^j+O(\Vert x,y\Vert^{r}),\\
			\dot{y}=-nx^{2n-1}+n\mu x^{n-1}y+\displaystyle\sum_{i+nj\geqslant 2n}\tilde{b}_{ij}x^iy^j+O(\Vert x,y\Vert^{r}).\\
		\end{array}
	\end{equation}
\end{proposition}
\noindent\textbf{Proof: }Note that the origin of system \eqref{eq1nHomo} is monodromic by Theorem \ref{TeoAndreevRestrito}. Conversely, considering the functions defined in Theorem \ref{TeoAndreevRestrito}, the change of variables $x\to x$, $y\to y-F(x)$, transforms system \eqref{eqPlanar} into
\begin{equation}\label{eqFNAndreev}
	\begin{array}{lcr}
		\dot{x}=y+y\tilde{P}(x,y),\\
		\dot{y}=f(x)+y\Phi(x)+y^2\tilde{Q}(x,y).\\
	\end{array}
\end{equation}
By Theorem \ref{TeoAndreevRestrito}, since the origin is monodromic, the above system can be rewritten as
\begin{equation*}
	\begin{array}{lcr}
		\dot{x}=y+\sum_{k\geqslant n+1}P_k(x,y)+O(\Vert x,y\Vert^{r}),\\
		\bigskip
		\dot{y}=ax^{2n-1}+\tilde{b}x^{n-1}y+\sum_{k\geqslant 2n}Q_k(x,y)+O(\Vert x,y\Vert^{r}),\\
	\end{array}
\end{equation*}
where $P_k(x,y),Q_k(x,y)\in\mathcal{P}^{(1,n)}_{k}$, $a\neq 0$, $\tilde{b}\in\mathbb{R}$ and $2n+1<r$. Note that the origin has Andreev number $n$. By Lemma \ref{LemaMonodromiaDelta}, we have $\Delta=\tilde{b}^2+4an<0$. Define $A=\frac{\Delta}{4n^2}$. The mapping $$\Psi(x,y)=\left(|A|^{\frac{1}{2(n-1)}}x, |A|^{\frac{1}{2(n-1)}}\left(y-\frac{\tilde{b}}{2n}x^{n}\right)\right),$$
is a diffeomorphism, and it transforms \eqref{eqFNAndreev} into
\begin{equation*}
	\begin{array}{lcr}
		\dot{x}=y+\mu x^n+\sum_{k\geqslant n+1}\tilde{P}_k(x,y)+O(\Vert x,y\Vert^{r}),\\
		\bigskip
		\dot{y}=-nx^{2n-1}+n\mu x^{n-1}y+\sum_{k\geqslant 2n}\tilde{Q}_k(x,y)+O(\Vert x,y\Vert^{r}),\\
	\end{array}
\end{equation*}
where $\tilde{P}_k(x,y),\tilde{Q}_k(x,y)\in\mathcal{P}^{(1,n)}_{k}$ and $\mu=\frac{\tilde{b}}{2n|A|^{\frac{1}{2}}}$.\qed

\begin{remark}\label{obsbeta=n-1}
	The expression of $\mu$ given above yields the following implication: Consider system \eqref{eqPlanar} with a monodromic nilpotent singular point. It satisfies the monodromy condition $\beta=n-1$ (item (ii) from Theorem \ref{TeoAndreevRestrito}) if and only if $\mu\neq0$ in canonical form \eqref{eq1nHomo}. In other words, the conditions (i) and (ii) in Theorem \ref{TeoAndreevRestrito} are preserved in the canonical form \eqref{eq1nHomo}.
\end{remark}

Now we introduce the generalized polar coordinates. Given $n$ a positive integer, consider the following differential system
$$\begin{array}{lcr}
	\dfrac{du}{d\theta}=-v;\\
	\dfrac{dv}{d\theta}=u^{2n-1},
\end{array}$$
with initial conditions $u(0)=1$, $v(0)=0$. Let $\mbox{Cs}\,\theta$ and $\mbox{Sn}\,\theta$ be the solutions of this Cauchy problem. We shall call these functions the \emph{generalized trigonometric functions}. The following result present their most useful properties and its proof can be found in \cite{LyapunovP,GasullTorre,AlvarezGasull1}.

\begin{proposition}\label{PropoGenTrigoProp}
	The following holds:
	\begin{itemize}
		\item[a)] {\rm Cs}\,$\theta$ and {\rm Sn}\,$\theta$ are $T$-periodic functions, where 
		$$T=2\sqrt{\dfrac{\pi}{n}}\dfrac{\Gamma(\frac{1}{2n})}{\Gamma(\frac{n+1}{2n})};$$
		\item[b)] $\mbox{\rm Cs}^{2n}\theta+n\mbox{\rm Sn}^2\theta=1$;
		\item[c)] {\rm Cs}\,$\theta$ is an even function and {\rm Sn}\,$\theta$ is an odd function;
		\item[d)] {\rm Cs}\,$(\frac{T}{2}+\theta)=-\mbox{\rm Cs}\,\theta$ and {\rm Sn}\,$(\frac{T}{2}+\theta)=-\mbox{\rm Sn}\,\theta$;
		\item[e)] {\rm Cs}\,$(\frac{T}{2}-\theta)=-\mbox{\rm Cs}\,\theta$ and {\rm Sn}\,$(\frac{T}{2}-\theta)=\mbox{\rm Sn}\,\theta$;
		\item[f)] $\displaystyle\int_{0}^{\theta}\mbox{\rm Sn}\,\varphi\;\mbox{\rm Cs}^q\varphi \;d\varphi=-\dfrac{\mbox{\rm Cs}^{q+1}\theta}{q+1}$;
		\item[g)] $\displaystyle\int_{0}^{\theta}\mbox{\rm Sn}^p\varphi\;\mbox{\rm Cs}^{2n-1}\varphi\;d\varphi=\dfrac{\mbox{\rm Sn}^{p+1}\theta}{p+1}$;
		\item[h)] $\displaystyle\int_{0}^{T}\mbox{\rm Sn}^p\varphi\;\mbox{\rm Cs}^q\varphi\;d\varphi=0$, if either $p$ or $q$ are odd;
		\item[i)] $\displaystyle\int_{0}^{T}\mbox{\rm Sn}^p\varphi\;\mbox{\rm Cs}^q\varphi\;d\varphi=\frac{2}{\sqrt{n^{p+1}}}\dfrac{\Gamma\left(\frac{p+1}{2}\right)\Gamma\left(\frac{q+1}{2n}\right)}{\Gamma\left(\frac{p+1}{2}+\frac{q+1}{2n}\right)}$, if both $p$ and $q$ are even;
	\end{itemize}
\end{proposition}

When working with planar vector fields having a nilpotent singular point, a useful change of variables based on the above defined functions is the following:
$$x=\rho\,\mbox{Cs}\,\theta,\;y=\rho^n\,\mbox{Sn}\,\theta.$$
We refer to this change of variables as the \emph{generalized polar coordinates}.

We now turn our attention to systems in family \eqref{eq1nHomo}. Applying the generalized polar coordinate change yields:
\begin{equation*}
	\begin{array}{lcr}
		\dot{\rho}=\rho^{n}\left(\mu\mbox{Cs}^{n-1}\theta-(n-1)\mbox{Cs}^{2n-1}\theta\,\mbox{Sn}\,\theta+P(\rho,\theta)+O(\rho^{r})\right),\\
		\dot{\theta}=-n\rho^{n-1}\left(1-(n-1)\mbox{Sn}^2\theta+Q(\rho,\theta)+O(\rho^{r})\right),
	\end{array}
\end{equation*}
where $P(\rho,\theta)$ and $Q(\rho,\theta)$ are given by:
\begin{eqnarray}
	P(\rho,\theta)&=&\displaystyle\sum_{i+nj\geqslant n+1}\tilde{a}_{ij}\rho^{i+nj-n}\mbox{Cs}^{i+2n-1}\theta\,\mbox{Sn}^j\theta\nonumber\\
	&&+\displaystyle\sum_{i+nj\geqslant 2n}\tilde{b}_{ij}\rho^{i+nj-2n+1}\mbox{Cs}^{i}\theta\,\mbox{Sn}^{j+1}\theta,\nonumber\\
	Q(\rho,\theta)&=&\displaystyle\sum_{i+nj\geqslant n+1}\tilde{a}_{ij}\rho^{i+nj-n}\mbox{Cs}^{i}\theta\,\mbox{Sn}^{j+1}\theta\nonumber\\
	&&-\displaystyle\sum_{i+nj\geqslant 2n}\frac{\tilde{b}_{ij}}{n}\rho^{i+nj-2n+1}\mbox{Cs}^{i+1}\theta\,\mbox{Sn}^{j}\theta.\nonumber
\end{eqnarray}
We then obtain the following differential equation:
\begin{equation}\label{eqRgen}
	\dfrac{d\rho}{d\theta}=R(\rho,\theta)=-\rho\dfrac{\left(\mu\mbox{Cs}^{n-1}\theta-(n-1)\mbox{Cs}^{2n-1}\theta\,\mbox{Sn}\,\theta+P(\rho,\theta)+O(\rho^{r})\right)}{n\left(1-(n-1)\mbox{Sn}^2\theta+Q(\rho,\theta)+O(\rho^{r})\right)}.
\end{equation}
Using Proposition \ref{PropoGenTrigoProp} item (b), $1-(n-1)\mbox{Sn}^2\theta=\mbox{Cs}^{2n}\theta+\mbox{Sn}^2\theta$, and thus the denominator of \eqref{eqRgen} does not vanish at $\rho=0$. Thus $R(\rho,\theta)$ is a $C^r$ function in a neighborhood of $\rho=0$. Let $\tilde{\rho}(\theta,\rho_0)$ be the solution of \eqref{eqRgen} with initial condition $\tilde{\rho}(0,\rho_0)=\rho_0$. We can expand it as a Taylor polynomial, i.e.
$$\tilde{\rho}(\theta,\rho_0)=\displaystyle\sum_{k=1}^{r}v_k(\theta)\rho_0^k+O(\rho^{r}),$$
and define the so-called \emph{focal values}.

\begin{definition}
	Consider function $\tilde{\rho}(\theta,\rho_0)$ associated to equation \eqref{eqRgen} defined above. The \emph{displacement function} is defined by $$d(\rho_0)=\tilde{\rho}(T,\rho_0)-\rho_0.$$
	The \emph{focal values} are the coefficients of the Taylor expansion of $d(\rho_0)$ given by $v_1(T)-1$ and $v_k(T)$ for $k\geqslant 2$.
\end{definition}

For the origin to be a center for system \eqref{eq1nHomo} all the focal values must be null. In the literature, for the non-degenerate singular point, it is well known that the first non-zero focal value is the coefficient of an odd power of $\rho_0$ (see \cite[Lemma 5, page 243]{AndronovB} or \cite[Proposition 3.1.4, page 94]{Romanovski}). In \cite{LiuLi3Ord}, the authors studied system \eqref{eq1nHomo} for $n=2$ and concluded that the first non-zero focal value is the coefficient of an even power of $\rho_0$. As a matter of fact, a more general result is true.

\begin{proposition}\label{PropoPrimeiroVFnnulogen}
	For system \eqref{eq1nHomo} with monodromic nilpotent singular points having Andreev number $n$, the first non-zero focal value is the coefficient of a power of $\rho_0$ whose parity is the same as the parity of $n$.
\end{proposition}

A proof of this proposition for the analytic case can be found in \cite[pages 50 to 54]{LyapunovP}. The $C^r$ case also holds, and the key argument to prove Proposition \ref{PropoPrimeiroVFnnulogen} is that the differential equation \eqref{eqRgen} is invariant under the following transformations:
$$\rho\mapsto -\rho,\;\theta\mapsto\theta+\frac{T}{2},\;\mbox{for $n$ odd},$$
$$\rho\mapsto -\rho,\;\theta\mapsto -\left(\theta+\frac{T}{2}\right),\;\mbox{for $n$ even}.$$
It is sufficient to follow the steps of \cite[Section 24 of Chapter IX]{AndronovB} to prove the proposition. 

\begin{lemma}
	The coefficient $v_1$ in function $\tilde{\rho}(\theta,\rho_0)$ associated to equation \eqref{eqRgen} is given by:
	\begin{equation}\label{eqv1gen}
		\exp\left(-\displaystyle\int_{0}^{\theta}\dfrac{\left(\mu\mbox{Cs}^{n-1}\varphi-(n-1)\mbox{Cs}^{2n-1}\varphi\,\mbox{Sn}\,\varphi\right)}{n\left(1-(n-1)\mbox{Sn}^2\varphi\right)}d\varphi\right).
	\end{equation}
\end{lemma}
\noindent\textbf{Proof: }Substituting the Taylor expansion of $\tilde{\rho}(\theta,\rho_0)$ in the equation \eqref{eqRgen} and comparing the coefficients of $\rho_0$ yields the result.\qed

\begin{lemma}\label{Lemav1gen}
	Consider function $\tilde{\rho}(\theta,\rho_0)$ associated to equation \eqref{eqRgen}. Then $v_1(T)=1$ for $n$ even and $v_1(T)=e^{-\mu A^2}$ for $n$ odd, where $A\neq 0$.
\end{lemma}
\noindent\textbf{Proof: }Using the properties described in Proposition \ref{PropoGenTrigoProp}, we compute the expression \eqref{eqv1gen} for $\theta=T$:
$$v_1(T)=\exp\left(-\mu\displaystyle\int_{0}^{T}\dfrac{\left(\mbox{Cs}^{n-1}\theta\right)}{n\left(1-(n-1)\mbox{Sn}^2\theta\right)}d\theta\right).$$

Let $G(\theta)$ be the integrand of the above expression. For $n$ odd, $G(\theta)$ is a non-negative even function. Thus, the above integral is not zero, and then $v_1(T)=e^{-\mu A^2}$ for some $A\neq 0$. If $n$ is even, $G(\theta)$ is an odd $T$-periodic function. Therefore $v_1(T)=1$.\qed

\begin{theorem}\label{TeoFocoNilpotenteForte}
	Suppose that the $C^r$ system \eqref{eqPlanar} has nilpotent monodromic singular point with odd Andreev number $n$. If the system satisfies the monodromy condition $\beta=n-1$ then the origin is a nilpotent focus.
\end{theorem}
\noindent\textbf{Proof: }Consider system \eqref{eqPlanar} under the hypothesis of Theorem \ref{TeoFocoNilpotenteForte}. By Proposition \ref{PropoPreformMono}, it can be written as system \eqref{eq1nHomo}. By Remark \ref{obsbeta=n-1}, we have $\mu\neq0$ in canonical form \eqref{eq1nHomo}. Lemma \ref{Lemav1gen} implies that $v_1(T)\neq 1$ for $\mu\neq 0$. Therefore, the origin must be a nilpotent focus. \qed
\medskip

Although the Andreev number is invariant by orbital equivalence (see \cite{GarciaFiI} for a proof), the number $\beta$ is not. The following system \cite[Example 13]{GarciaGineInte},
\begin{equation*}
	\begin{array}{lcr}
		\dot{x}=y+a_1x^9+a_2x^6y-a_1(54a_1^2-6a_2+b_2)x^3y^2+a_4y^3,\\
		\dot{y}=-x^{11}-9a_1x^8y+b_2x^5y^2+3a_1(18a_1^2+b_2)x^2y^3,\\
	\end{array}
\end{equation*}
has Andreev number $n=6$ and $\beta=14$. Performing an analytic change of variables, we can transform the above system into the form \eqref{eqFNAndreev}:
\begin{equation*}
	\begin{array}{lcr}
		\dot{x}=y+O(\Vert x,y\Vert^{13}),\\
		\dot{y}=-x^{11}(1+O(x))+yx^{11}+O(\Vert x,y\Vert^{13}),\\
	\end{array}
\end{equation*}
for which the Andreev number is $n=6$ and $\beta=11$. However, note that for both systems, the monodromy condition $\beta>n-1$, i.e. condition (ii) in Theorem \ref{TeoAndreevRestrito}, is satisfied. In fact, using the canonical form \eqref{eqFNAndreev}, it is possible to prove that the monodromy conditions (i) and (ii) in Theorem \ref{TeoAndreevRestrito} are invariant by local diffeomorphisms. Before stating the result, we need the following definition of quasi-homogeneous vector fields.

\begin{definition}
	A planar map $F:\mathbb{R}^2\to\mathbb{R}^2$, $F(x,y)=(F_1(x,y),F_2(x,y))$ is \emph{$(t_1,t_2)$-quasi-homogeneous of weighted degree $k$} if $F_1(x,y)\in\mathcal{P}^{(t_1,t_2)}_{k+t_1}$ and $F_2(x,y)\in\mathcal{P}^{(t_1,t_2)}_{k+t_2}$. The vector space of all $(t_1,t_2)$-quasi-\-homogeneous maps of weighted degree $k$ is denoted by $\mathcal{Q}^{(t_1,t_2)}_{k}$. A planar vector field $X$ associated to a system $\mathbf{\dot{x}}=F(\mathbf{x})$ is \emph{$(t_1,t_2)$-quasi-\-homogeneous of weighted degree $k$} if $F\in\mathcal{Q}^{(t_1,t_2)}_{k}$. In this case, we write $X\in\mathcal{Q}^{(t_1,t_2)}_{k}$
\end{definition}

One can argue that a more natural definition would be taking $F_1,F_2\in\mathcal{P}^{(t_1,t_2)}_{k}$. However, the above definition is better suited for our study and is also used in several papers, for instance \cite{AlgabaInt, GarciaFiI}.

\begin{remark}
	The vector spaces $\mathcal{P}^{(t_1,t_2)}_{k}$ and $\mathcal{Q}^{(t_1,t_2)}_{k}$ are well-behaved. Some of their properties that are useful in the understanding of the next proof are listed below:
	\begin{itemize}
		\item If $f\in \mathcal{P}^{(t_1,t_2)}_{k}$ then $\dfrac{\partial f}{\partial x}\in \mathcal{P}^{(t_1,t_2)}_{k-t_1}$ and $\dfrac{\partial f}{\partial y}\in \mathcal{P}^{(t_1,t_2)}_{k-t_2}$;
		\item If $f\in \mathcal{P}^{(t_1,t_2)}_{k_1}$ and $g\in \mathcal{P}^{(t_1,t_2)}_{k_2}$ then $fg\in \mathcal{P}^{(t_1,t_2)}_{k_1+k_2}$;
		\item If $\varphi\in \mathcal{Q}^{(t_1,t_2)}_{k_1}$ and $X\in \mathcal{Q}^{(t_1,t_2)}_{k_2}$ then $d\varphi\cdot X\in \mathcal{Q}^{(t_1,t_2)}_{k_1+k_2}$.
	\end{itemize}
\end{remark}

\begin{proposition}\label{Propobeta=n-1invariante}
	For $C^r$ system \eqref{eqPlanar} having an isolated monodromic nilpotent singular point at the origin, with Andreev number $n$ such that $2n-1<r$, the monodromy conditions (i) and (ii) in Theorem \ref{TeoAndreevRestrito} are invariant by local diffeomorphisms.
\end{proposition}
\noindent\textbf{Proof: }Let $X$ and $Y$ be vector fields having an isolated nilpotent singular point satisfying the monodromy conditions (i) and (ii) in Theorem \ref{TeoAndreevRestrito}, respectively. Suppose that there exists a local diffeomorphism $(\bar{x},\bar{y})=\varphi(x,y)$ that transforms $X$ into $Y$. Without loss of generality, we can assume that the associated systems are in canonical form \eqref{eqFNAndreev}. We have that 
\begin{equation}\label{eqConjugation1}
	d\varphi\cdot Y=X\circ\varphi.
\end{equation}
We write $\varphi(x,y)=(\sum_{j+k\geqslant 1}a_{jk}x^jy^k,\sum_{j+k\geqslant 1}b_{jk}x^jy^k)+O(\Vert x,y\Vert^{r})$. The homogeneous terms of degree 1 in both sides of \eqref{eqConjugation1} must satisfy the following equation:
$$\left(\begin{array}{l}
	a_{10}y\\
	b_{10}y
\end{array}\right)=\left(\begin{array}{c}
	b_{10}x+b_{01}y\\
	0
\end{array}\right).$$
Thus, $a_{10}=b_{01}$ and $b_{10}=0$. Note that $a_{10}$ must not be zero or else $\varphi$ is not a diffeomorphism in a neighborhood of the origin. 

Now, we write $X=X_{n-1}+\bar{X}$, $Y=Y_{n-1}+\bar{Y}$ where $X,Y\in\mathcal{Q}_{n-1}^{(1,n)}$ and $\bar{X},\bar{Y}$ are $C^r$ vector fields containing the $(1,n)$-quasi-homogeneous terms of degree greater than $n-1$. We also write $\varphi=\sum_{k=-(n-1)}^{n-1}\varphi_{k}+\tilde{\varphi}$ where $\varphi_{k}\in\mathcal{Q}_{k}^{(1,n)}$ and $\tilde{\varphi}$ is a $C^r$ map containing the $(1,n)$-quasi-homogeneous terms of degree greater than $n-1$. Note that, for $k=1,\dots,n-1$, $\varphi_{-k}$ is given by:
$$\varphi_{-k}(x,y)=(0,b_{n-k,0}x^{n-k})\;\mbox{and}\;\varphi_0=(a_{10}x,a_{10}y+b_{n0}x^n).$$
Since $X$ and $Y$ satisfy the monodromy conditions (i) and (ii) in Theorem \ref{TeoAndreevRestrito}, respectively, we have
$$X_{n-1}=\left(\begin{array}{l}
	y\\
	ax^{2n-1}
\end{array}\right),\;Y_{n-1}=\left(\begin{array}{l}
	y\\
	Ax^{2n-1}+B x^{n-1}y
\end{array}\right),$$
for some $a,A<0$ and $B\neq 0$. 
Let $k\in\{1,\dots,n-1\}$ be such that $\varphi_{-k}$ is the first non-zero $(1,n)$-quasi-homogeneous term in $\varphi$, i.e. $\varphi=\varphi_{-k}+\dots+\tilde{\varphi}$. Then the lowest $(1,n)$-quasi-homogeneous term in the expression $d\varphi\cdot Y$ is $d\varphi_{-k}\cdot Y_{n-1}$, whose weighted degree is $n-k-1$, i.e.
$$d\varphi_{-k}\cdot Y_{n-1}=\left(\begin{array}{l}
	0\\
	(n-k)b_{n-k,0}x^{n-k-1}y
\end{array}\right).$$
Now, the $(1,n)$-quasi-homogeneous term of weighted degree $n-k-1$ in $X\circ\varphi$ is $(b_{n-k}x^{n-k},0)^T$. Equating this expression to the one above, we obtain that $b_{n-k,0}=0$, which means that $\varphi_{-k}\equiv 0$. But this is a contradiction. Thus, $\varphi=\sum_{k=0}^{n-1}\varphi_{k}+\tilde{\varphi}$.

We go back to \eqref{eqConjugation1} once more to complete the proof. The lowest $(1,n)$-quasi-homogeneous term in the expression $d\varphi\cdot Y$ is $d\varphi_{0}\cdot Y_{n-1}$, whose weighted degree is $n-1$. Denoting by $(X\circ\varphi)_{n-1}$ the $(1,n)$-quasi-homogeneous term of degree $n-1$ in $X\circ\varphi$, it follows that $d\varphi_0\cdot Y_{n-1}=(X\circ\varphi)_{n-1}$, that is
$$\left(\begin{array}{l}
	a_{10}y\\
	Aa_{10}x^{2n-1}+(nb_{n0}+a_{10}B)x^{n-1}y
\end{array}\right)=\left(\begin{array}{l}
	a_{10}y+b_{n0}x^n\\
	aa_{10}^{2n-1}x^{2n-1}
\end{array}\right).$$
For this equation to be satisfied, we must have $b_{n0}=0$ and $a_{10}B=0$. Since $a_{10}\neq 0\neq B$, this is not possible. We conclude that if $\varphi$ is a local diffeomorphism, it must preserve the monodromy conditions.
\qed

\begin{remark}
	The above result holds for formal settings, that is, formal vector fields and formal changes of variables $\varphi$ such that $d\varphi(0)\neq 0$. 
\end{remark}

\section{Auxiliary Linear Operators}\label{SecLinearOperators}

In this section, we prove auxiliary results regarding some useful linear operators that show up frequently when working with vector fields having nilpotent singular points. These particular linear operators provide some interesting results on normal forms and methods involving formal series.

We denote by $H^{(n)}_{(3,1)}$ the vector space of homogeneous polynomials of degree $n$ in three variables. A basis for this vector space is given by the collection of all monomials $x^jy^kz^l$ for which $j,k,l\geqslant 0$ and $j+k+l=n$. Given a polynomial $p\in H^{(n)}_{(3,1)}$, let $\langle p\rangle$ be the vector subspace spanned by $p$. We consider the following linear operator:
$$\begingroup
\setlength\arraycolsep{0pt}
T_n\colon\begin{array}[t]{c >{{}}c<{{}} c}
	H^{(n)}_{(3,1)} & \to & H^{(n)}_{(3,1)} \\
	\noalign{\medskip} p & \mapsto & y\dfrac{\partial p}{\partial x}-\lambda z\dfrac{\partial p}{\partial z},
\end{array}
\endgroup$$

where $\lambda\neq 0$.
\begin{lemma}\label{LemaKerTn}
	The kernel of linear operator $T_n$ is given by $\ker T_n=\langle y^n\rangle$.
\end{lemma}
\noindent\textbf{Demo:} Clearly $T_n(y^n)=0$. Let $p=\sum_{j+k+l=n}p_{j,k,l}x^jy^kz^l$. If $p\in\ker T_n$, then $T_n(p)=0$, i.e.:
$$\sum_{j+k+l=n}((j+1)p_{j+1,k-1,l}-\lambda lp_{j,k,l})x^jy^kz^l=0.$$
Consider $l\neq 0$. For $k=0$, the coefficient of $x^jz^l$ in the left-hand side of this equation is $-\lambda lp_{j,0,l}$ which implies that $p_{j,0,l}=0$. For $k=1$, the equation yields $p_{j,1,l}=0$. Proceeding step by step until index $k=n-1$, we conclude that $p_{j,k,l}=0$. Now, for $l=0$, the coefficient of monomial $x^jy^k$ in the left-hand side expression above is $(j+1)p_{j+1,k-1,0}$, which implies that $p_{j,k,0}=0$ for $j>0$. Thus $p=p_{0,n,0}y^n\in\langle y^n\rangle$.
\qed

\begin{lemma}\label{LemaTp+q}
	For every $q\in H^{(n)}_{(3,1)}$, there is a choice of $p\in H^{(n)}_{(3,1)}$ such that $T_n(p)+q\in\langle x^n\rangle$.
\end{lemma}
\noindent\textbf{Proof: }It is sufficient to prove $H^{(n)}_{(3,1)}=\langle x^n\rangle\oplus \text{Im} T_n$. By Lemma \ref{LemaKerTn}, the codimension of Im$T_n$ is 1. Since $\langle x^n\rangle$ has dimension 1, we only need to prove that the intersection of both vector subspaces is $\{0\}$. We have $T_n(p)\in\langle x^n\rangle$ if and only if $y\dfrac{\partial p}{\partial x}-\lambda z\dfrac{\partial p}{\partial z}=\alpha x^n$, which is only possible when $\dfrac{\partial p}{\partial x}=\dfrac{\partial p}{\partial z}=\alpha=0$. Therefore $T_n(p)\in\langle x^n\rangle$ if and only if $p\in\langle y^n\rangle=\ker T_n$. The result follows. \qed

\medskip
Now, we turn our attentions to the following linear operator:
$$\begingroup
\setlength\arraycolsep{0pt}
L_n\colon\begin{array}[t]{c >{{}}c<{{}} c}
	H^{(n)}_{(3,1)} & \to & H^{(n)}_{(3,1)} \\
	\noalign{\medskip} p & \mapsto & y\dfrac{\partial p}{\partial x}-\lambda z\dfrac{\partial p}{\partial z}+\lambda p,
\end{array}
\endgroup$$
where $\lambda\neq 0$.

\begin{lemma}\label{LemaKerLnJacobi}
	The kernel of linear operator $L_n$ is given by $\ker L_n=\langle y^{n-1}z\rangle$.
\end{lemma}
\noindent\textbf{Demo:} Trivially $L_n(y^{n-1}z)=0$. Let $p=\sum_{j+k+l=n}p_{j,k,l}x^jy^kz^l$. If $p\in\ker L_n$, then $L_n(p)=0$ which equates to:
$$\sum_{j+k+l=n}((j+1)p_{j+1,k-1,l}-\lambda (l-1)p_{j,k,l})x^jy^kz^l=0.$$
Consider $l\neq 1$. For $k=0$, the coefficient of $x^jz^l$ in the left-hand side of this equation is ${-\lambda (l-1)p_{j,0,l}}$ which implies that $p_{j,0,l}=0$. For $k=1$, the equation yields $p_{j,1,l}=0$. Proceeding step by step until index $k=n-1$, we conclude that $p_{j,k,l}=0$. Now, for $l=1$, the coefficient of monomial $x^jy^kz$ in the left-hand side expression above is $(j+1)p_{j+1,k-1,1}$, which implies that $p_{j,k,1}=0$ for $j>0$. Thus $p=p_{0,n-1,1}y^{n-1}z\in\langle y^{n-1}z\rangle$.\qed

\begin{lemma}
	For every $q\in H^{(n)}_{(3,1)}$, there is a choice of $p\in H^{(n)}_{(3,1)}$ such that $L_n(p)+q\in\langle x^{n-1}z\rangle$.
\end{lemma}
\noindent\textbf{Demo: }Again, we prove that $H^{(n)}_{(3,1)}=\langle x^{n-1}z\rangle\oplus \text{Im} L_n$. Lemma \ref{LemaKerLnJacobi} tells us that the codimension of Im$L_n$ is 1. As in the proof of Lemma \ref{LemaTp+q}, it is enough to prove that the intersection of $\langle x^{n-1}z\rangle$ and Im$L_n$ is $\{0\}$. We have that $L_n(p)\in\langle x^{n-1}z\rangle$ if and only if $y\dfrac{\partial p}{\partial x}-\lambda z\dfrac{\partial p}{\partial z}+\lambda p=\alpha x^{n-1}z$. The coefficient of  $x^{n-1}z$ in the left-hand side of this equation is zero, therefore $\alpha=0$. Hence $L_n(p)\in\langle x^{n-1}z\rangle$ if and only if $p\in\langle y^{n-1}z\rangle=\ker L_n$. The result holds. \qed

\begin{remark}\label{ObsOperator}
	Note that if we consider similar operators
	$$\tilde{T}_n:H^{(n)}_{(3,1)}\to H^{(n)}_{(3,1)},\; p\mapsto x\dfrac{\partial p}{\partial y}-\lambda z\dfrac{\partial p}{\partial z}, \mbox{and }$$
	$$\tilde{L}_n:H^{(n)}_{(3,1)}\to H^{(n)}_{(3,1)},\; p\mapsto y\dfrac{\partial p}{\partial x}-\lambda z\dfrac{\partial p}{\partial z}+\lambda p,$$
	just by exchanging $x$ and $y$ in the expressions above yields analogous results. More precisely, we have $\ker\tilde{T}_n=\langle x^n\rangle$, $\ker \tilde{L}_n=\langle x^{n-1}z\rangle$  and $H^{(n)}_{(3,1)}=\langle y^n\rangle\oplus\mbox{\rm Im}\tilde{T}_n=\langle y^{n-1}z\rangle\oplus\mbox{\rm Im}\tilde{L}_n$.
\end{remark}

\section{Normal Forms}\label{SecNormalForms}


An important step in the study of vector fields having nilpotent singular points is to find normal forms for its associated differential system. We first look at the Normal Form theory to search the most practical normal forms for system \eqref{eq1}. The following theorem, which is proven in section 1.5 of \cite{Zhitomirskii} gives us a starting point:

\begin{theorem}[Zhitomirskii Normal Form]\label{TeoFormaNormalZhi} Consider a formal system of differential equations of form
	$$\dot{\mathbf{x}}=A\mathbf{x}+X(\mathbf{x}),$$
	where $\mathbf{x}\in\mathbb{R}^n$, $A\in M_n(\mathbb{R})$ and $j^1X(0)=0$. Then there exists a germ at $0$ of diffeomorphism $\mathbf{x}\to\mathbf{x}+\varphi(\mathbf{x})$, $j^1\varphi(0)=0$, that transforms the above system into the following:
	$$\dot{\mathbf{x}}=A\mathbf{x}+f(\mathbf{x}),$$
	where $j^1f(0)=0$ and 	
	\begin{equation}\label{condFN}
		A^tf(\mathbf{x})-df(\mathbf{x})A^t\mathbf{x}=0.
	\end{equation}
\end{theorem}

We apply Theorem \ref{TeoFormaNormalZhi} to system \eqref{eq1} and obtain its correspondent normal form. System \eqref{eq1} is written in the desired form $\dot{\mathbf{x}}=A\mathbf{x}+X(\mathbf{x})$ where ${\mathbf{x}\in\mathbb{R}^3}$ with $A=$\begin{tiny}
	$\left(\begin{array}{ccc}
		0 & 1 & 0\\
		0 & 0 & 0\\
		0 & 0 & -\lambda
	\end{array}\right)$\end{tiny}. Writting function $f(\mathbf{x})$ in Theorem \ref{TeoFormaNormalZhi} as $f(\mathbf{x})=(f_1(\mathbf{x}),f_2(\mathbf{x}),f_3(\mathbf{x}))^t$, we can solve the equation \eqref{condFN} to obtain its expression. Since
$$A^tf(\mathbf{x})=\left(\begin{array}{c}
	0\\
	f_1(\mathbf{x})\\
	-\lambda f_3(\mathbf{x})
\end{array}\right),\;\mbox{and}\;df(\mathbf{x})A^t\mathbf{x}=\left(\begin{array}{c}
	x\dfrac{\partial f_1}{\partial y}(\mathbf{x})-\lambda z \dfrac{\partial f_1}{\partial z}(\mathbf{x})\\
	x\dfrac{\partial f_2}{\partial y}(\mathbf{x})-\lambda z \dfrac{\partial f_2}{\partial z}(\mathbf{x})\\
	x\dfrac{\partial f_3}{\partial y}(\mathbf{x})-\lambda z \dfrac{\partial f_3}{\partial z}(\mathbf{x})
\end{array}\right).$$
equation \eqref{condFN} becomes
\begin{eqnarray}\label{condFN1}
	-x\dfrac{\partial f_1}{\partial y}(\mathbf{x})+\lambda z \dfrac{\partial f_1}{\partial z}(\mathbf{x})=0,\nonumber\\
	f_1(\mathbf{x})-x\dfrac{\partial f_2}{\partial y}(\mathbf{x})+\lambda z\dfrac{\partial f_2}{\partial z}(\mathbf{x})=0,\\
	-\lambda f_3(\mathbf{x})-x\dfrac{\partial f_3}{\partial y}(\mathbf{x})+\lambda z\dfrac{\partial f_3}{\partial z}(\mathbf{x})=0.\nonumber
\end{eqnarray}
Let us write $f_i(\mathbf{x})=\sum_{n\geqslant 2}f^{n}_i(\mathbf{x})$ where $f^{n}_i(\mathbf{x})$ are homogeneous polynomials of degree $n$ for $i=1,2,3$. The first equation in \eqref{condFN1} implies
$$x\dfrac{\partial f^{n}_1}{\partial y}(\mathbf{x})-\lambda z \dfrac{\partial f^{n}_1}{\partial z}(\mathbf{x})=0,\;\mbox{for }n\geqslant 2.$$
Recalling Remark \ref{ObsOperator} and Lemma \ref{LemaKerTn}, by the above equation we obtain that for $n\geqslant 2$, $f^{n}_1$ must be in the kernel of linear operator $\tilde{T}_n$, i. e. $f^{n}_1\in\langle x^n\rangle$. Thus, we conclude that $f_1(\mathbf{x})=xP_1(x)$ with $P_1(0)=0$.

Now, we turn our attention to the third equation in \eqref{condFN1}. It is equivalent to the following equation:
$$x\dfrac{\partial f^{n}_3}{\partial y}(\mathbf{x})-\lambda z\dfrac{\partial f^{n}_3}{\partial z}(\mathbf{x})+\lambda f^{n}_3(\mathbf{x})=0,\;\mbox{for }n\geqslant 2,$$
which, by Lemma \ref{LemaKerLnJacobi} and Remark \ref{ObsOperator}, is equivallent to
$$f^{n}_3\in\ker\tilde{L}_n=\langle x^{n-1}z\rangle,$$
for $n\geqslant 2$. Hence, $f_3(\mathbf{x})=zR_1(x)$ where $R_1(0)=0$. It remains to find a more precise expression for function $f_2(\mathbf{x})$. The second equation in \eqref{condFN1} implies
$$x\dfrac{\partial f^{n}_2}{\partial y}(\mathbf{x})-\lambda z\dfrac{\partial f^{n}_2}{\partial z}(\mathbf{x})=f^{n}_1(\mathbf{x}),\;\mbox{for }n\geqslant 2.$$
To continue this investigation, it is convenient to write the functions $f_1(\mathbf{x}),f_2(\mathbf{x})$ as power series expansions: $f_i(\mathbf{x})=\sum_{j+k+l\geqslant 2}^{}F_{j,k,l}^{i}x^jy^kz^l$ for $i=1,2$. To simplify the arguments, we adopt the convention $F^{i}_{j,k,l}=0$ if any of the indexes $j,k,l$ is negative. Thus, the above equation is equivalent to
$$\sum_{j+k+l\geqslant 2}F_{j,k,l}^{1}x^jy^kz^{l}-\sum_{j+k+l\geqslant 2}^{}kF_{j,k,l}^{2}x^{j+1}y^{k-1}z^l+\sum_{j+k+l\geqslant 2}^{}\lambda lF_{j,k,l}^{2}x^jy^kz^{l}=0,$$
which yields:
\begin{equation}\label{condFN2} F^{1}_{j,k,l}-(k+1)F_{j-1,k+1,l}^{2}+\lambda lF_{j,k,l}^{2}=0.
\end{equation}
We know that $F^{1}_{j,k,l}=0$ for $k^2+l^2\neq 0$. For $l\neq 0$, equation \eqref{condFN2} implies that all coeficients $F^{2}_{j,k,l}$ are null. For $l=0$, \eqref{condFN2} becomes
$$F^{1}_{j,k,0}-(k+1)F_{j-1,k+1,0}^{2}=0.$$
Thus, the only coefficients of $f_2(\mathbf{x})$ not necessarily zero are $F^{2}_{j,1,0}$ and $F^{2}_{j,0,0}$. In fact, by the above equation, the coefficients $F^{2}_{j,1,0}$ must satisfy $F_{j,1,0}^{2}=F^{1}_{j+1,0,0}$. Moreover, the coefficients $F^{2}_{j,0,0}$ have no restriction presented by equation \eqref{condFN2}. Therefore, $f_2(\mathbf{x})=Q_2(x)+yP_1(x)$ where $j^1Q_2(0)=0$. Hence, we have proven the following result.

\begin{theorem}[Nilpotent Normal Form in $\mathbb{R}^3$]
	For system \eqref{eq1} having a nilpotent singular point at the origin, there exist a formal change of variables that transforms it into the formal normal form
	\begin{equation}\label{eqN1}
		\begin{array}{lr}
			\dot{x}=y+xP_1(x),\\
			\dot{y}=Q_2(x)+yP_1(x),\\
			\dot{z}=-\lambda z+zR_1(x).
		\end{array}
	\end{equation}
	for which $P_1(0)=j^1Q_2(0)=R_1(0)=0$.
\end{theorem}

However, we were not able to prove that the normal form \eqref{eqN1} is analytic, i.e. that the series $P_1,Q_2$ and $R_1$ are convergent. That does not mean that this normal form is not useful. We remark that the above normal form has $z=0$ as an invariant surface which is a center manifold and the first two components are decoupled from the third.

\begin{remark}
	The Zhitomirskii normal form for planar systems \eqref{eqPlanar} having a nilpotent singular point is
	\begin{equation}\label{eqN1Planar}
		\begin{array}{lr}
			\dot{x}=y+xP_1(x),\\
			\dot{y}=Q_2(x)+yP_1(x),
		\end{array}
	\end{equation}
	where $j^1Q_2(0)=P_1(0)=0$ (see Theorem 1.8.6, page 37 in \cite{Zhitomirskii}), which is very similar to its three-dimensional counterpart. From this form it is possible to arrive at the Liénard Canonical Form 
	\begin{equation}\label{eqLienardPlanar}
		\begin{array}{lr}
			\dot{x}=-y,\\
			\dot{y}=Q_2(x)+yP_1(x),
		\end{array}
	\end{equation}
	where $j^1Q_2(0)=P_1(0)=0$. In \cite{Zoladek} it is proven that this form is actually analytic. From \eqref{eqLienardPlanar}, it is possible to solve the planar nilpotent center problem, since the origin of \eqref{eqLienardPlanar} is a center if and only if $P_1(x)$ is an odd function \cite{BerthierMoussu}. Furthermore the authors of \cite{Zoladek} proved that the Liénard Canonical Form \eqref{eqLienardPlanar} is integrable if and only if $P_1(x)\equiv 0$.
\end{remark}

\section{Nilpotent centers and First Integrals}\label{SecFI}

It is well-known that the center problem for systems having a Hopf singularity is equivalent to the existence of a first integral at the point \cite{Bibikov, Adam, Queiroz}. In this section we search for some link between integrability and nilpotent centers on center manifolds. 

\begin{definition}
	Let $X$ be a vector field and $H:U\subset\mathbb{R}^k\to\mathbb{R}$ be a $C^1$-function not locally constant. If $H$ satisfies
	$$XH\equiv 0,$$
	then $H$ is a \emph{first integral} for $X$ and in this case $X$ is \emph{integrable}.
\end{definition}

In the remainder of this section, we use the following notation: given a function $F:\mathbb{R}^3\to\mathbb{R}$, we denote by $F_n$ the homogeneous part of degree $n$ in its power series expansion.

\begin{theorem}\label{PropoXH=xn}
	Consider a vector field $X$ associated to system \eqref{eq1} having a nilpotent singular point. Then there exists a formal series $H(x,y,z)$ $=y^2+\sum_{n\geqslant 3}H_n(x,y,z)$ such that ${XH=\sum_{n\geqslant 4}\omega_nx^n}$.
\end{theorem}
\noindent\textbf{Proof: }The Lie derivative $XH$ is given by:
\begin{small}
	\begin{eqnarray}
		XH &=&\left(y+P_2+P_3+\dots\right)\left(\dfrac{\partial H_3}{\partial x}+\dfrac{\partial H_4}{\partial x}+\dots\right)\nonumber\\
		&&+\left(-\lambda z+R_2+R_3+\dots\right)\left(\dfrac{\partial H_3}{\partial z}+\dfrac{\partial H_4}{\partial z}+\dots\right)\nonumber\\
		&&+\left(Q_2+Q_3+\dots\right)\left(2y+\dfrac{\partial H_3}{\partial y}+\dfrac{\partial H_4}{\partial y}+\dots\right).\nonumber
	\end{eqnarray}
\end{small}

Note that $y\dfrac{\partial H_n}{\partial x}-\lambda z\dfrac{\partial H_n}{\partial z}$ is a homogeneous polynomial of degree $n$. Hence, rewriting the above expression by organizing the homogeneous terms of degree $n$ in brackets, we have that
\begin{small}
	\begin{eqnarray}
		XH &=&\left(y\dfrac{\partial H_3}{\partial x}-\lambda z\dfrac{\partial H_3}{\partial z}+2yQ_2\right)+\left(y\dfrac{\partial H_4}{\partial x}-\lambda z\dfrac{\partial H_4}{\partial z}+F_4\right)\nonumber\\
		&&+\left(y\dfrac{\partial H_5}{\partial x}-\lambda z\dfrac{\partial H_5}{\partial z}+F_5\right)+\dots=\nonumber\\
		&=&\sum_{n\geqslant 3}\left(y\dfrac{\partial H_n}{\partial x}-\lambda z\dfrac{\partial H_n}{\partial z}+F_n\right)=\sum_{n\geqslant 3}T_n(H_n)+F_n,\nonumber
	\end{eqnarray}
\end{small}
where $F_n\in H^{(n)}_{(3,1)}$ are obtained by the homogeneous parts of $P,Q,R$ and $H$ of degree less than $n$. By Lemma \ref{LemaTp+q}, we can choose $H_n$ such that $T_n(H_n)+F_n=\omega_nx^n$, $\omega_n\in\mathbb{R}$. Moreover, for $n=3$, $F_3=2yQ_2\notin\langle x^3\rangle$ which means that there is a choice of $H_3$ such that $T_3(H_3)=-2yQ_2$. Then, by making the suitable choices, $XH=\sum_{n\geqslant 4}\omega_nx^n$.\qed
\medskip

Considering the normal form \eqref{eqN1}, the following result gives a characterization of the shape of a formal first integral.

\begin{lemma}\label{LemaHndependez}
	If $H$ is a formal first integral for the normal form \eqref{eqN1}, then $\dfrac{\partial H}{\partial z}\equiv 0$, that is $H=H(x,y)$.
\end{lemma}
\noindent\textbf{Proof: }Let $H(x,y,z)=\sum_{j+k+l\geqslant 1}p_{jkl}x^jy^kz^l$ be a formal first integral for \eqref{eqN1}. We assume $p_{jkl}=0$ if any of the indexes $j,k,l$ is negative. The following equation holds:
$$\dfrac{\partial H}{\partial x}(y+xP_1(x))+\dfrac{\partial H}{\partial y}\left(Q_2(x)+yP_1(x)\right)+\dfrac{\partial H}{\partial z}(-\lambda z+zR_1(x))=0.$$
The coefficients of the monomials $z^l$ in the left-hand side of the above equation are given by $-\lambda lp_{0,0,l}$. Hence, $p_{0,0,l}=0$ for $l\in\mathbb{N}$.

Now, we study the coefficients of monomials $x^jy^kz^l$ with $l\neq 0$ in the left-hand side of the above equation. We conclude that these coefficients satisfy the following relation:
\begin{eqnarray}\label{eqdemo1}
	(j+1)p_{j+1,k-1,l}+\sum_{m=2}^{j}(j-m+1)p_{j-m+1,k,l}a_{m-1}&\nonumber\\
	+\sum_{m=2}^{j}(k+1)p_{j-m,k+1,l}b_m&\nonumber\\
	+\sum_{m=1}^{j}kp_{j-m,k,l}a_{m}-\lambda lp_{j,k,l}+\sum_{m=1}^{j+k-1}lp_{j-m,k,l}c_{m}&=0,
\end{eqnarray}
where $a_m,b_m,c_m$ are the coefficients of $x^m$ in the expression of $P_1,Q_2,R_1$ respectively. Note that the summations above are expressed in terms of the coefficients $p_{j_1,j_2,l}$ with $j_1+j_2<j+k$. We now prove that $p_{jkl}=0$ if $l\neq 0$ by induction over $j+k$:

Suppose $j+k=1$. Since $p_{j_1,j_2,l}=0$ for $j_1+j_2<1$, equation \eqref{eqdemo1} becomes 
\begin{equation}\label{eqdemo2}
	(j+1)p_{j+1,k-1,l}-\lambda lp_{j,k,l}=0.
\end{equation}
We conclude that $p_{0,1,l}=p_{1,0,l}=0$. Now suppose that for $j+k<n$, $p_{j,k,l}=0$. Then, for $j+k=n$, equation $\eqref{eqdemo1}$ becomes \eqref{eqdemo2}. For $(j,k,l)=(n,0,l)$, we have by \eqref{eqdemo2} that $p_{n,0,l}=0$. Hence, equation \eqref{eqdemo2} sets the following chain of implications:
$$p_{n,0,l}=0\Rightarrow p_{n-1,1,l}=0\Rightarrow p_{n-2,2,l}=0\Rightarrow\dots\Rightarrow p_{0,n,l}=0.$$
Therefore $p_{j,k,l}=0$ for $l\neq 0$ and the result follows.\qed
\medskip

The previous lemma lets us conclude that formal integrability of the normal form \eqref{eqN1} is essentially formal integrability of the planar normal form \eqref{eqN1Planar}. Hence, we proceed to study the integrability of the formal system \eqref{eqN1Planar}. 

One of the consequences of Lemma \ref{LemaHndependez} is that if system \eqref{eq1} is formally integrable, then it has a formal first integral $H$ such that $j^2H(0)=y^2$. In fact, since \eqref{eq1} is transformed into \eqref{eqN1} via near-identity changes of variables and the formal integrability of system \eqref{eq1} is equivalent to the formal integrability of system \eqref{eqN1}, it is enough to verify this statement for system \eqref{eqN1Planar}. But this follows from the proof of Theorem 1 in paper \cite{Chavarriga}. Therefore, the quantities $\omega_n$ in Theorem \ref{PropoXH=xn} present obstructions for the system \eqref{eq1} to be analytically or formally integrable. For integrable systems, the monodromic singular point must be a center, since in this case the restricted system is also integrable. However not all nilpotent centers are formally integrable. For instance, consider the following system:
\begin{equation}\label{exNotIntegrableCenter}
	\begin{array}{lr}
		\dot{x}=y+x^2,\\
		\dot{y}=-x^3,\\
		\dot{z}=-\lambda z.
	\end{array}
\end{equation}
which has $z=0$ as a center manifold and its restriction is a time-reversible system which implies that the origin is a center on the center manifold. If we try to construct a formal first integral $H(x,y,z)$ for system we obtatin $\omega_5=2$. Note that the obstruction is a coefficient of odd power of $x$ in $XH$. If it was an even power, we would not have a nilpotent center, as the next result holds.

\begin{theorem}\label{PropoXHNFoco}
	Let $X$ be the vector field associated to system \eqref{eq1} having a nilpotent singular point and $H$ be a formal series as in Theorem \ref{PropoXH=xn}. If there exists $n\in\mathbb{N}$ such that $j^{2n}XH(0)=\omega_{2n}x^{2n}$ with $\omega_{2n}\neq 0$, then the origin cannot be a center on the center manifold.
\end{theorem}

\noindent\textbf{Proof: }Let $H(x,y)=y^2+\sum_{j=3}^{\infty}H_j(x,y)$ be the formal series given by Theorem \ref{PropoXH=xn}. Consider its $2n$-jet, i. e.: $\tilde{H}(x,y)=y^2+\sum_{j=3}^{2n}H_j(x,y)$. Thus, $\tilde{H}$ is an analytic function such that $j^{2n}X\tilde{H}(0)=\omega_{2n}x^{2n}$.

Let $z=h(x,y)$ be a $C^{2n+3}$-parametrization of a center manifold for $X$. For the restricted system, we have:
\begin{equation}\label{condPropoXHNFoco}
	X\vert_{z=h(x,y)}\tilde{H}\vert_{z=h(x,y)}=X\tilde{H}\vert_{z=h(x,y)}=\omega_{2n}x^{2n}+r(x,y),
\end{equation}
where $j^{2n}r(0)=0$. Denote by $X_{W^c}$ and $\tilde{H}_{W^c}$ the restriction of $X$ and $\tilde{H}$ to the center manifold $z=h(x,y)$. 

Without loss of generality suppose $\omega_{2n}<0$. The same conclusions will follow for $\omega_{2n}>0$ considering vector field $-X$.

Since $j^2\tilde{H}_{W^c}(0)=y^2$ and $j^{2n}X_{W^c}\tilde{H}_{W^c}(0)=\omega_{2n}x^{2n}$, there is a neighborhood $U$ of the origin on the center manifold such that $\tilde{H}_{W^c}\geqslant 0$ and $X_{W^c}\tilde{H}_{W^c}\leqslant 0$. Let $p\in U$ and $\gamma_p(t)$ the trajectory of the restriction system \eqref{eq1} to $W^c$ with initial point $\gamma_p(0)=p$. Note that $\tilde{H}_{W^c}$ is non-increasing along $\gamma_p$. If $\gamma_p$ is periodic and wholly contained in $U$, there is $T\in\mathbb{R}$ such that $\gamma_p(0)=\gamma_p(T)$ and thus $\tilde{H}_{W^c}\circ\gamma_p(0)=\tilde{H}_{W^c}\circ\gamma_p(T)$. By continuity of $\tilde{H}_{W^c}$, we have $X_{W^c}\tilde{H}_{W^c}\vert_{\gamma_p}\equiv 0$. Therefore, any closed orbits of $X_{W^c}$ inside $U$ are contained in the level sets of $\tilde{H}_{W^c}$. Thus, if the origin is a center, there is a neighborhood of it for which $X_{W^c}\tilde{H}_{W^c}$ is null. But this contradicts \eqref{condPropoXHNFoco}.\qed

\medskip

By Proposition \ref{PropoMono2}, the example \eqref{exNotIntegrableCenter} has a monodromic singular point with Andreev number $n=2$ and satisfies monodromy condition $\beta=n-1$ in Theorem \ref{TeoAndreevRestrito}. The lack of formal integrability for this system is not an exception. More precisely, the monodromy condition $\beta=n-1$ is an obstruction for formal integrability as we will see in the next results.

\begin{lemma}\label{Lemabeta=n-1notIntFN}
	Consider system \eqref{eqN1Planar} having a monodromic singular point satisfying monodromy condition $\beta=n-1$ in Theorem \ref{TeoAndreevRestrito}. Then it is not formally integrable.
\end{lemma}
\noindent\textbf{Proof: }For system \eqref{eqN1Planar}, the functions $F(x),f(x)$ and $\Phi(x)$ in Theorem \ref{TeoAndreevRestrito} are given by 
$$F(x)=-xP_1(x),\;f(x)=Q_2(x)-xP_1^2(x),\;\Phi(x)=2P_1(x)+x\dfrac{\partial P_1(x)}{\partial x}.$$
Writing $P_1(x)=\sum_{l\geqslant 1}a_lx^{l}$ and $Q_2(x)=\sum_{l\geqslant 2}b_lx^l$, we have:
$$f(x)=b_2x^2+\sum_{l\geqslant 3}(b_l-A_{l-1})x^{l},\;\mbox{ and }\;\Phi(x)=\sum_{l\geqslant 1}(2+l)a_lx^{l},$$
where $A_l=\sum_{i=1}^{l-1}a_ia_{l-i}$. Since monodromy condition $\beta=n-1$ in Theorem \ref{TeoAndreevRestrito} is satisfied, there exist a positive integer $n\geqslant 2$ such that ${(b_{2n-1}-A_{2n-2})<0}$ is the first non-zero coefficient in the power series expansion of $f(x)$ and $a_{n-1}$ is the first non-zero coefficient in the power series expansion of $P_1(x)$. Note that $A_k=0$ for $k<2(n-1)$, $b_k=0$ for $k<2n-1$ and ${A_{2n-2}=a_{n-1}^2}$ which implies $b_{2n-1}-a_{n-1}^2<0$.

Suppose there exists $H(x,y)$ a formal first integral for system \eqref{eqN1Planar}. Thus, we can write $H(x,y)=y^2+\sum_{j+k\geqslant 3}p_{jk}x^jy^k$. Then, the following equation holds:
\begin{equation}\label{21}
	y\dfrac{\partial H}{\partial x}+Q_2(x)\dfrac{\partial H}{\partial y}+P_1(x)\left(x\dfrac{\partial H}{\partial x}+y\dfrac{\partial H}{\partial y}\right)\equiv 0.
\end{equation}
Studying the above equation by expanding its left-hand side in a power series, we conclude that the coefficients of the monomials $x^jy^k$ are given by the left-hand side of the following equation:
\begin{equation}\label{22}
	(j+1)p_{j+1,k-1}+\sum_{l=0}^{j-2}\left((k+1)b_{j-l}p_{l,k+1}+a_{j-l}(l+k)p_{l,k}\right)=0.
\end{equation}
The above equation implies that for $0<j\leqslant n-1$, $p_{j,k}=0$. Now, substituting $j=n-1$ in \eqref{22} yields $np_{n,k-1}+ka_{n-1}p_{0,k}=0$, and so $$p_{n,k-1}=-\dfrac{ka_{n-1}p_{0,k}}{n}.$$

Substituting $k=1$ in \eqref{22} yields:
\begin{equation}\label{24}
	(j+1)p_{j+1,0}+\sum_{l=0}^{j-2}\left(2b_{j-l}p_{l,2}+a_{j-l}(l+1)p_{l,1}\right)=0,
\end{equation}
and for $j\leqslant n$, it follows that $p_{j,0}=0$. Studying the above equation for $n\leqslant j\leqslant 2n-1$, we conclude that $p_{j,0}=0$, for $j<2n$, and $2np_{2n,0}+2b_{2n-1}+a_{n-1}(n+1)p_{n,1}=0$,
i.e.
$$p_{2n,0}=-\dfrac{\left(2b_{2n-1}+a_{n-1}(n+1)p_{n,1}\right)}{2n}.$$

Now, we can explicitly compute the coefficient of $x^{3n-1}$ in the left-hand side of equation \eqref{21} denoted by $\omega_{3n-1}$:
\begin{eqnarray}
\omega_{3n-1}&=&\sum_{j=2}^{3n-3}(b_{3n-1-j}p_{j,1}+ja_{3n-1-j}p_{j,0})=\nonumber\\
&=&\sum_{j=n}^{2n}(b_{3n-1-j}p_{j,1}+ja_{3n-1-j}p_{j,0})=\nonumber\\
&=&b_{2n-1}p_{n,1}+2na_{n-1}p_{2n,0}=-\dfrac{2a_{n-1}(n+1)(a_{n-1}^2-b_{2n-1})}{n}.\nonumber
\end{eqnarray}
Formal integrability implies that $\omega_{3n-1}=0$, but this is not possible since $a_{n-1}\neq 0$ and $b_{2n-1}-a_{n-1}^2<0$. Therefore, the result holds.\qed

\begin{remark}
	Using the notation from Theorem \ref{TeoAndreevRestrito}, the quantity $\omega_{3n-1}$ given in the previous proof is given by:
	$$\omega_{3n-1}=-\dfrac{2ab}{n}.$$
\end{remark}

\begin{proposition}
	Consider system \eqref{eqPlanar} having a monodromic singular point satisfying monodromy condition $\beta=n-1$ in Theorem \ref{TeoAndreevRestrito}. Then it cannot admit formal first integral.
\end{proposition}
\noindent\textbf{Proof: }By Proposition \ref{Propobeta=n-1invariante}, the monodromy condition $\beta=n-1$ is invariant under formal near-identity changes of variables. Moreover, the formal integrability of system \eqref{eqPlanar} would imply the formal integrability of its normal form \eqref{eqN1Planar}. Hence, by Lemma \ref{Lemabeta=n-1notIntFN}, system \eqref{eqPlanar} cannot be formally integrable.\qed

\begin{theorem}\label{Teobeta=n-1notInt}
	Consider system \eqref{eq1} having a monodromic singular point such that its restriction to a center manifold satisfies monodromy condition $\beta=n-1$ in Theorem \ref{TeoAndreevRestrito}. Then it cannot admit formal first integral.
\end{theorem}
\noindent\textbf{Proof: }The result follows from the fact that the formal integrability of system \eqref{eq1} implies the formal integrability of the associated normal form \eqref{eqN1} and, by Lemma \ref{LemaHndependez}, the formal integrability of the bidimensional normal form \eqref{eqN1Planar}. By Lemma \ref{Lemabeta=n-1notIntFN}, system \eqref{eq1} cannot be formally integrable.\qed

\medskip

The proof of Lemma \ref{Lemabeta=n-1notIntFN} has an interesting implication that helps us in the study of the Nilpotent Center Problem for three-dimensional system \eqref{eq1}. More precisely, we have the following result.

\begin{theorem}\label{TeoAndreevimparbeta=n-1NotCenter}
	Suppose that the origin of system \eqref{eq1} is monodromic with odd Andreev number $n$ and satisfies the monodromy condition $\beta=n-1$. Then the origin cannot be a center on a center manifold.
\end{theorem}
\noindent\textbf{Proof: }Let $X$ be the vector field associated to system \eqref{eq1}, $Y$ the formal vector field associated to the normal form \eqref{eqN1} and $\varphi$ the formal near-identity change of variables that transforms \eqref{eq1} into \eqref{eqN1}, i.e. $Y\circ\varphi=d\varphi\cdot X$. Let us denote $(u,v,w)=\varphi(x,y,z)$. By the proof of Lemma \ref{Lemabeta=n-1notIntFN}, it is possible to construct a formal series $\tilde{H}(u,v,w)$ with $j^2\tilde{H}(0)=v^2$ such that $Y\tilde{H}=\omega_{3n-1}u^{3n-1}+O(u^{3n})$. Let $H=\tilde{H}\circ\varphi$. We have
$$XH=\langle\nabla H,X\rangle=\langle\nabla\tilde{H}\cdot d\varphi,d\varphi^{-1}\cdot Y\circ\varphi\rangle=$$
$$=\langle\nabla\tilde{H},Y\rangle\circ\varphi=\omega_{3n-1}x^{3n-1}+O(x^{3n}).$$ Since $n$ is odd, $3n-1$ is an even number. Thus, by Theorem \ref{PropoXHNFoco}, the origin cannot be a center on a center manifold.\qed

\medskip

Note that an alternative proof of Theorem \ref{TeoAndreevimparbeta=n-1NotCenter} can be given using Theorem \ref{TeoFocoNilpotenteForte}. The monodromy condition $\beta>n-1$ also interferes with the integrability of the normal forms \eqref{eqN1} and \eqref{eqN1Planar}. More precisely, the next results classify the normal form for the integrable nilpotent singular points in $\mathbb{R}^3$.

\begin{theorem}\label{TeoFNIntegravel}
	Consider system \eqref{eqN1Planar} having a monodromic singular point. If it admits formal first integral $H(x,y)$, then either $P_1(x)\equiv 0$ or $m=2sn-1$ for some $s\in\mathbb{N}$.
\end{theorem}
\noindent\textbf{Proof: }For system \eqref{eqN1Planar}, the functions $F(x),f(x)$ and $\Phi(x)$ in Theorem \ref{TeoAndreevRestrito} are given by 
$$F(x)=-xP_1(x),\;f(x)=Q_2(x)-xP_1^2(x),\;\Phi(x)=2P_1(x)+x\dfrac{\partial P_1(x)}{\partial x}.$$
Writing $P_1(x)=\sum_{l\geqslant 1}a_lx^{l}$ and $Q_2(x)=\sum_{l\geqslant 2}b_lx^l$, we have:
$$f(x)=b_2x^2+\sum_{l\geqslant 3}(b_l-A_{l-1})x^{l},\;\mbox{ and }\;\Phi(x)=\sum_{l\geqslant 1}(2+l)a_lx^{l},$$
where $A_l=\sum_{i=1}^{l-1}a_ia_{l-i}$. The monodromy conditions for system \eqref{eqN1Planar} imply that $b_2=0$ and the first nonzero coefficient in the power series expansion of $f(x)$ is $(b_{2n-1}-A_{2n-2})<0$. Furthermore, either $P_1(x)\equiv 0$, $m>n-1$, or $m=n-1$ where $m$ is the index of the first nonzero coefficient $a_m$ of $P_1(x)$. By Lemma \ref{Lemabeta=n-1notIntFN}, integrability can only occur in the cases $m>n-1$ or $P_1(x)\equiv 0$. 

For $P_1(x)\equiv 0$, $H(x,y)=y^2-\sum_{l\geqslant 2}\frac{2 b_l}{l+1}x^{l+1}$ is a formal first integral for system \eqref{eqN1Planar}. Now, we suppose that $m>n-1$. Then $A_j=0$ for $j<2m$, $b_j=0$ for $j<2n-1$ and $b_{2n-1}<0$. If there exists $H(x,y)$ a formal first integral for system \eqref{eqN1Planar}, we can write $H(x,y)=y^2+\sum_{j+k\geqslant 3}p_{jk}x^jy^k$ and the equation \eqref{21} holds. Moreover the coefficients $p_{j,k}$ satisfy equations \eqref{22} and \eqref{24}.

Suppose $n\leqslant m<2n-1$. Studying equations \eqref{22} and \eqref{24}, we conclude that $p_{j,k}=0$ for $0<j\leqslant m$ and $p_{j,0}=0$ for $j\leqslant m+1$. Substituting $j=m$ in \eqref{22} yields $p_{m+1,k-1}=-\dfrac{ka_{m}p_{0,k}}{m+1}$. By \eqref{24} for $m+1\leqslant j\leqslant 2n-1$, we have that:
$$p_{j,0,}=0,\;\mbox{for}\;j<2n,\;p_{m+2,0}=-\dfrac{2b_{m+1}}{m+2},\;\mbox{and}\;p_{2n,0}=-\dfrac{2b_{2n-1}}{2n}.$$
Note that $p_{m+2,0}\neq 0$ only if $m=2n-2$. With this information, we compute the coefficient $\omega_{2n+m}$ of $x^{2n+m}$ in $XH$, which is given by:
\begin{eqnarray}
\omega_{2n+m}&=&\sum_{j=2}^{2n+m-2}(b_{2n+m-j}p_{j,1}+ja_{2n+m-j}p_{j,0})=\nonumber\\
&=&\sum_{j=m+1}^{2n}(b_{2n+m-j}p_{j,1}+ja_{2n+m-j}p_{j,0}),\nonumber
\end{eqnarray} 
which simplifies to:
$$\omega_{2n+m}=2na_{m}p_{2n,0}+b_{2n-1}p_{m+1,1}=-\dfrac{2a_mb_{2n-1}(m+2)}{m+1}.$$
Since $a_m\neq 0$ and $b_{2n-1}<0$, it follows that $\omega_{2m+n}\neq 0$ and therefore, system \eqref{eqN1Planar} cannot be integrable.

Now, suppose $m>2n-1$ and $m\neq 2sn-1$ for $s\in\mathbb{N}$. Thus, there exists $r\in\mathbb{N}$ such that $2rn-1<m<2(r+1)n-1$. By \eqref{22} we have that $p_{j,k}=0$ for $0<j\leqslant 2n-1$. 
For $2n-1\leqslant j<m$, equation \eqref{22} becomes:
\begin{equation}\label{25}
	(j+1)p_{j+1,k-1}+\sum_{l=0}^{j-2n+1}(k+1)b_{j-l}p_{l,k+1}=0.
\end{equation}

Moreover, the quantities $\omega_{2sn-1}$ for $s\leqslant r+1$ are given by:
\begin{equation}\label{26}
	\omega_{2sn-1}=\sum_{j=2n}^{2(s-1)n}(b_{2sn-1-j}p_{j,1}).
\end{equation}
Substituting $k=2$ in \eqref{25}, yields:
\begin{equation}\label{27}
	(j+1)p_{j+1,1}+\sum_{l=0}^{j-2n+1}\left(3b_{j-l}p_{l,3}\right)=0.
\end{equation}
For $j=2n-1$ we have $p_{2n,1}=-\dfrac{3b_{2n-1}p_{0,3}}{2n}$ and substituting this expression in \eqref{26} for $s=2$, we conclude that $\omega_{4n-1}=-\dfrac{3b_{2n-1}^2p_{0,3}}{2n}$. Formal integrability implies $p_{0,3}=0$. As a consequence, $p_{j,1}=0$ for $j\leqslant 2n$.

We turn again to equation \eqref{24} with $j\leqslant 2n-1, p_{j,0}=0$ and $p_{2n,0}=-\dfrac{b_{2n-1}}{n}$. The quantity $\omega_{2n+m}$ is given by:

\begin{eqnarray}
\omega_{2n+m}&=&\sum_{j=2}^{2n+m-2}(b_{2n+m-j}p_{j,1}+ja_{2n+m-j}p_{j,0})=\nonumber\\
&=&-2a_{m}b_{2n-1}+\sum_{j=2n+1}^{m+1}(b_{2n+m-j}p_{j,1}).\nonumber
\end{eqnarray}

Our goal is to prove that $\omega_{2n+m}$ cannot be zero. Recalling that $p_{j,1}=0$ for $j\leqslant 2n$ and $p_{j,k}=0$ for $0<j<2n$. By \eqref{27}, for $j\leqslant 4n-2$, we have:

$$(j+1)p_{j+1,1}+\sum_{l=0}^{j-2n+1}3b_{j-l}p_{l,3}=(j+1)p_{j+1,1}+3b_{j}p_{0,3}=0,$$

and thus $p_{j,1}=0$ for $j\leqslant 4n-1$. Moreover, substituting $j=4n-1$ in \eqref{27}, yields
$$4np_{4n,1}+\sum_{l=0}^{2n}3b_{4n-1-l}p_{l,3}=4np_{4n,1}+3b_{2n-1}p_{2n,3}=0,$$
and from \eqref{25},
$$2np_{2n,3}+5b_{2n-1}p_{0,5}=0.$$

Hence, we obtain $\omega_{6n-1}=b_{2n-1}p_{4n,1}$ and the formal integrability implies that $p_{4n,1}=0$ and $p_{0,5}=0$.

Suppose, by induction, that for $j\leqslant 2sn$, $p_{j,1}=0$ and $p_{0,3}=p_{0,5}=\dots=p_{0,2s+1}=0$, for $s<r$. Using equations \eqref{27} and \eqref{25} repeatedly, we obtain for $j\leqslant 2(s+1)n-2$,

$$(j+1)p_{j+1,1}+\sum_{l=0}^{j-2n+1}3b_{j-l}p_{l,3}=$$
$$=(j+1)p_{j+1,1}+\sum_{l=0}^{j-2n+1}3b_{j-l}\left(-\dfrac{1}{l}\sum_{t=0}^{l-2n}5b_{l-1-t}p_{t,5}\right).$$

The expressions of the coefficients $p_{j,k}$, by \eqref{25}, are given in terms of nested summations whose indexes vary in intervals of natural numbers that are each instance $2n$ units less than the previous one until eventually fit in the interval $0\leqslant j\leqslant 2n-1$. Then, all coefficients, $p_{j,k}$ are null for $j>0$. Thus, the coefficient $p_{j+1,1}$ is a linear combination of the coefficients $p_{0,3},p_{0,5},\cdots,p_{0,2s+1}$. By the induction hypothesis, we conclude that $p_{j,1}=0$ for $j\leqslant 2(s+1)n-1$.

Computing $p_{2(s+1)n,1}$ and $\omega_{2(s+1)n-1}$, we conclude that the system \eqref{eqN1} can only be formally integrable if $p_{0,2(s+1)+1}=0$ and consequentially $p_{2(s+1)n,1}=0$.

Thus, using \eqref{27}, we compute $p_{m+1,1}$:
\begin{equation*}
	(m+1)p_{m+1,1}+2a_m+\sum_{l=1}^{m-2}\left(3b_{m-l}p_{l,3}\right)=0,
\end{equation*}
which implies that $p_{m+1,1}=-\dfrac{2a_m}{m+1}$. Then $\omega_{2n+m}$ becomes:
$$\omega_{2n+m}=-2a_{m}b_{2n-1}-\dfrac{2a_mb_{2n-1}}{m+1},$$
which is not zero. Hence, system \eqref{eqN1} with monodromic singular point can only be integrable if $P_1(x)\equiv 0$ or $m=2sn-1$. \qed

\medskip

Theorem \eqref{TeoFNIntegravel} gives us a formal normal form for integrable systems \eqref{eq1}. The same result holds for planar systems and it was proven by a different method in \cite[Theorem 4]{GarciaGineInte}. Although it is clear that the case $P_1(x)\equiv 0$ is always formally integrable, there are few examples of the case $m=2sn-1$. One of those is the following system:
\begin{equation*}
	\begin{array}{lr}
		\dot{x}=y+\dfrac{x^4}{5},\\
		\dot{y}=-x^3+\dfrac{x^{7}}{25}+\dfrac{x^3y}{5},\\
		\dot{z}=-\lambda z.
	\end{array}
\end{equation*}
Here $n=2$ and $m=3=2n-1$. The above system has $H(x,y)=\left(y+\frac{x^4}{4}-1\right)\exp\left(y-\frac{x^4}{20}\right)$ as a first integral. Note that $j^2H(0)=-1+\frac{y^2}{2}$.

Even though we have considered formal integrability, the following result, proven by Mattei and Moussu in \cite{MatteiMoussu}, tells us that formal integrability and analytic integrability are the same.

\begin{theorem}[Theorem A in \cite{MatteiMoussu}]
	Let $F:\mathbb{R}^n\to\mathbb{R}^n$ be an analytic vector field, $0$ an isolated singular point of $F$ and $H$ a formal first integral of $F$ such that $H$ is not a power of a formal series. Then, there exists a formal function $f:\mathbb{R}\to\mathbb{R}$, with $f(0)=0,f'(0)\neq 0$ such that $f\circ H$ is an analytic first integral of $F$.
\end{theorem}

\begin{remark}
	The condition that $H$ is not a power of a formal series, means that, any decomposition $H=h_1^{k_1}h_2^{k_2}\dots h_l^{k_l}$ in irreducible factors is such that the greatest common divisor of $k_1,\dots,k_l$ is $1$.
\end{remark}

\section{Applications}\label{SecApplications}
\subsection{Illustrative example}

We solve the Nilpotent Center Problem for the following system
\begin{equation}\label{eqKuklesQuad}
	\begin{array}{lcr}
		\dot{x}=y,\\
		\dot{y}=b_{101}xz+b_{{020}}y^2+b_{{011}}yz+b_{002}z^2,\\
		\dot{z}=-z+{x}^{2}+2xy+c_{020}y^2,
	\end{array}
\end{equation}
using a method based on the above results. We assume that $b_{101}<0$ so that, by Proposition \ref{PropoMono2}, the origin is monodromic with Andreev number $2$. Using Theorem \ref{PropoXH=xn}, we compute the quantities $\omega_n$ which are obstructions for the system to be integrable. We have that the first non-zero is
$$\omega_6=-\dfrac{2\,b_{011}b_{101}}{3}.$$
Thus, by Theorem \ref{PropoXHNFoco}, for the origin of \eqref{eqKuklesQuad} to be a center, we must have $b_{011}=0$. Under this condition, we compute the next non-zero $\omega_n$ which is
$$\omega_8=\dfrac{4\,c_{020}b_{101}^3}{5}.$$
Therefore, $b_{011}=c_{020}=0$ are necessary conditions for the origin of \eqref{eqKuklesQuad} to be a center on a center manifold. These conditions are in fact sufficient, since under those, \eqref{eqKuklesQuad} admits $z=x^2$ as a center manifold and the restricted system is
\begin{equation}
	\begin{array}{lcr}
		\dot{x}=y,\\
		\dot{y}=b_{101}x^3+b_{{020}}y^2+b_{002}x^4,\\
	\end{array}
\end{equation}
which is reversible. Hence the origin of \eqref{eqKuklesQuad} is a nilpotent center on a center manifold if and only if $b_{011}=c_{020}=0$.

\subsection{The Generalized Lorenz system}

The Generalized Lorenz system is one of the most studied three-dimensional systems in the literature since its dynamics are very rich. Among the particular cases of the Generalized Lorenz systems are the Lü and Chen systems. Its expression is given by 
\begin{equation}\label{eqLorenz}
	\begin{array}{lcr}
		\dot{x}=a(y-x),\\
		\dot{y}=bx+cy-xz,\\
		\dot{z}=dz+xy.
	\end{array}
\end{equation}
We consider $ad\neq 0$ for the singular points of system \eqref{eqLorenz} to be isolated. The origin is always a singular point whose type of singularity is determined by the parameters $a,b,c,d$. The nondegenerate case is studied in \cite{GarciaHopf,Queiroz} and the references therein. We investigate conditions on the parameters $a,b,c,d$ for which system \eqref{eqLorenz} has a nilpotent singular point.

The singular points of the above system are the origin, and $Q_{\pm}={(\pm\sqrt{-d(b+c)},\pm\sqrt{-d(b+c)},b+c)}$ for $d(b+c)<0$. The determinant of the Jacobian matrix of \eqref{eqLorenz} at point $Q_{\pm}$ is given by $2ad(b+c)$. Thus, $Q_{\pm}$ are nilpotent singular points only when $b+c=0$ which means they coincide with the origin. 

Now, the Jacobian matrix of \eqref{eqLorenz} at the origin is
$$\begin{small}
	\left(\begin{array}{rcr}
		-a & a & 0 \\
		b & c & 0 \\
		0 & 0 & d
	\end{array}
	\right),
\end{small}$$
and for the origin to be a nilpotent singular point, we must have $b+c=0$ and $c=a$. By means of the coordinate change $\bar{x}=y$, $\bar{y}=a(y-x)$, $\bar{z}=z$, dropping the bars, system \eqref{eqLorenz} becomes
\begin{equation}\label{eqLorenzNill}
	\begin{array}{lcr}
		\dot{x}=y-xz+\tfrac{1}{a}yz,\\
		\dot{y}=-axz+yz,\\
		\dot{z}=dz+x^2-\tfrac{1}{a}xy.
	\end{array}
\end{equation}
Theorem \ref{PropoXHNFoco} is powerful enough to solve the Nilpotent Center Problem for the above system. We compute the quantities $\omega_n$ (Theorem \ref{PropoXH=xn}) for system \eqref{eqLorenzNill} and obtain the first non-zero one:
$$\omega_6=-\frac {2\,a \left( 2\,a+d \right) }{3{d}^{3}}.$$
Thus, by Theorem \ref{PropoXHNFoco}, the origin can only be a nilpotent center on a center manifold if $d=-2a$. Under this condition, by Proposition \ref{PropoMono2}, the origin is monodromic with Andreev number $2$. Moreover, it satisfies monodromy condition $\beta>n-1$ from Theorem \ref{TeoAndreevRestrito}.

For $d=-2a$, the function $V(x,y,z)=x^2-\frac{2xy}{a}+\frac{y^2}{a^2}-2az$ defines an invariant surface $V\equiv 0$ for \eqref{eqLorenzNill} which is tangent to the $xy$-plane, thus, it is a center manifold for system
\eqref{eqLorenzNill}. Moreover, the restriction of the system to this center manifold is given by
\begin{equation*}
	\begin{array}{lcr}
		\dot{x}=y-xz+\dfrac{y}{2a^2}\left(x-\dfrac{y}{a}\right)^2,\\
		\\
		\dot{y}=-axz+\dfrac{y}{2a}\left(x-\dfrac{y}{a}\right)^2,
	\end{array}
\end{equation*}
which is a Hamiltonian system with Hamiltonian function $H(x,y)=y^2+\frac{x^4}{4}-{\frac {{x}^{3}y}{a}}+{\frac {3{x}^{2}{y}^{2}}{2{a}^{2}}}-{\frac {x{y}^{3}}{{a}^{3}}}+{\frac {{y}^{4}}{4{a}^{4}}}$. Thus, the origin is a nilpotent center on a center manifold. We conclude:

\begin{theorem}
	The origin of the Generalized Lorenz system \eqref{eqLorenz} is a nilpotent center on a center manifold if and only if $b=-a,\,c=a,\,d=-2a$.
\end{theorem}


\subsection{Hide-Skeldon-Acheson Dynamo system}

In paper \cite{HideSkeldon} the authors proposed a model for self-exciting dynamo action in which a Faraday disk and coil are arranged in series with either a capacitor or a motor. The proposed system for this model is the following
\begin{equation}\label{eqHideSkeldon}
	\begin{array}{lcr}
		\dot{x}=-x-\beta z+xy,\\
		\dot{y}=\alpha-\kappa y-\alpha x^2,\\
		\dot{y}=x-\lambda z.
	\end{array}
\end{equation}
The physical quantities that the parameters $\alpha,\beta,\lambda,\kappa$ represent are detailed in \cite{HideSkeldon, HideSkeldon2}. The above system has a singular point at $Q=(0,\frac{\alpha}{\kappa},0)$. The determinant of the Jacobian matrix of \eqref{eqHideSkeldon} at $Q$ is $ \left( \alpha-\kappa \right) \lambda-\beta\,\kappa$, and its characteristic polynomial has
$${\frac { \left( \lambda+1 \right) {\kappa}^{2}+ \left( \beta-\alpha+
		\lambda \right) \kappa-\lambda\,\alpha}{\kappa}}$$
as the coefficient of its linear term. Both these expressions must be zero for $Q$ to be a nilpotent singular point. This occurs for $\alpha=\kappa(\lambda+1)$ and $\beta=\lambda^2$. Now, under these conditions, we translate $Q$ to the origin and via the change of coordinates ${x=\lambda\bar{x}+\bar{y},\;y=\bar{z},\;z=\bar{x}}$, dropping the bars, we obtain the following system
\begin{equation}\label{eqHideSkeldon2}
	\begin{array}{lcr}
		\dot{x}=y,\\
		\dot{y}=\lambda xz+yz,\\
		\dot{y}=-\kappa z-\kappa (\lambda+1)(\lambda x+y)^2.
	\end{array}
\end{equation}
We now investigate whether the origin is monodromic or not. The inequality \eqref{ineqMono2} in Proposition \ref{PropoMono2} for system \eqref{eqHideSkeldon2} is given by:
$$-\lambda^3(\lambda+1)<0,$$
which is only satisfied for $\lambda<-1$ or $\lambda>0$. For $-1<\lambda<0$, by Theorem \ref{TeoAndreevRestrito} the origin of system \eqref{eqHideSkeldon2} is not monodromic (see the proof of Proposition \ref{PropoMono2}). Moreover, for $\lambda=-1$ or $\lambda=0$, \eqref{eqHideSkeldon2} does not have an isolated singular point since the $x$-axis is a line of equilibria in this case. Note that, by Proposition \ref{PropoMono2} we also know that when the origin of system \eqref{eqHideSkeldon2} is monodromic, it satisfies the condition $\beta>n-1$ in Theorem \ref{TeoAndreevRestrito}.

We compute the quantities $\omega_n$ described in Theorem \ref{PropoXH=xn}. The first non-zero is
$$\omega_6=-{\frac {2{\lambda}^{5} \left( \lambda+1 \right) ^{2} \left( -2\,		\lambda+3\,\kappa \right) }{3\kappa}},$$
and the origin can only be a nilpotent center when $\lambda=\frac{3\kappa}{2}$. Under this condition, the next non-zero quantity is $\omega_8$ which is also even indexed:
$$\omega_8=-{\frac {6561\,{\kappa}^{6} \left( 3\,\kappa+2 \right) ^{3}}{512}}.$$
We obtain that $\kappa=-\frac{2}{3}$ is the only value for which the origin could be a nilpotent center. However, that would imply that $\lambda=-1$ and thus the origin would not be an isolated singular point. In fact, for $\lambda=-1$, $z=0$ is a center manifold and the flow of \eqref{eqHideSkeldon2} is governed by the linear vector field $y\partial_x$ which does not have any periodic orbits. Thus, we have proven the following result.

\begin{theorem}
	The singular point $Q=(0,\frac{\alpha}{\kappa},0)$ is never a nilpotent center on a center manifold for system \eqref{eqHideSkeldon2}.
\end{theorem}

\section{Acknowledgments}

The first author is partially supported by S\~ao Paulo Research Foundation (FAPESP) grants 19/10269-3 and 18/19726-5. The second author is supported by S\~ao Paulo Research Foundation (FAPESP) grant 19/13040-7. 

\section*{Appendix}

\noindent\textbf{Proof of Theorem \ref{TeoAndreevRestrito}:  }Function $F$ is well defined by the Implicit Function Theorem. By making the change of variables $x=x,\; y=v+F(x)$, system \eqref{eqPlanarCr} becomes:
\begin{equation*}
	\begin{array}{lr}
		\dot{x}=v+F(x)+X_2(x,v+F(x))=v(1+X_1(x,v)),\\
		\dot{v}=\dot{y}-F'(x)\dot{x}=f(x)+\Phi(x)v+Y_0(x,v)v^2.
	\end{array}
\end{equation*}
We distinguish two cases.
\noindent Case $\alpha$ even: By means of linear parametrization, we can assume $a=1$. Consider the circunference $C_\delta=\{x^2+v^2=\delta^2\}$ with $\delta$ small enough for $X$ to be well defined on $C_\delta$. To compute the Poincaré index at the origin we consider the $v$-axis. The points $M_1=(\delta,0)$, $M_2=(-\delta,0)$, are such that $X(M_1)=(0,f(\delta))$ and $X(M_2)=(0,f(-\delta))$ are paralell to the $v$-axis, and are the only ones on $C_\delta$ with this property.

Since $\alpha$ is even, we have $f(\delta)f(-\delta)>0$. Evaluating the vector field near $M_1$ and $M_2$ following the curve $C_\delta$ counterclockwise, we see that the vector field near $M_1$ ceases to be paralell to the $v$-axis to be in the quadrant $x>0, v>0$. Meanwhile near $M_2$, the vector field enter the quadrant $x<0, v>0$. In other words, following $C_\delta$ counterclockwise, the vector field turns clockwise near $M_1$ and counterclockwise near $M_2$. Thus, the Poincaré index at the origin is $0$ and therefore it cannot be a monodromic singular point.

Case $\alpha$ odd: Considering $C_\delta$ once more, we apply the same argument as above to compute the Poincaré index at the origin. Now we have $f(\delta)f(-\delta)<0$. If $a>0$, then $f(\delta)>0$ and as we evaluate vector field $X$ along the curve $C_\delta$ counterclockwise it turns clockwise near $M_1$ and clockwise near $M_2$. Hence, Poincaré index at the origin is $-1$.


Now for $a<0$, then $f(\delta)<0$, and as we follow $C_\delta$ counterclockwise, the vector field turns counterclockwise near $M_1$ and counterclockwise near $M_2$. Therefore, Poincaré index at the origin is $1$. We conclude that monodromy can only occur for $a<0$.

Suppose $\alpha=2n-1$ and $a<0$. Making the change of variables $x=(-\frac{1}{a})^{\frac{1}{2n-2}}\bar{x},\; v=-(-\frac{1}{a})^{\frac{1}{2n-2}}\bar{y}$, and dropping the bars, yields the following system:

\begin{equation}\label{eqDemoTeoAndreev2}
	\begin{array}{lr}
		\dot{x}=y(-1+X_1(x,y)),\\
		\dot{y}=g(x)+\Phi(x)y+Y_0(x,y)y^2,
	\end{array}
\end{equation}

where $g(x)=x^{2n-1}+O(x^{2n})$. We have three cases to consider:

Case $\beta<n-1$: We now make the generalized polar coordinate change. System \eqref{eqDemoTeoAndreev2} becomes
\begin{equation*}
	\begin{array}{lr}
		\dot{\rho}=b\text{Cs}^\beta(\theta)\text{Sn}^2(\theta)\rho^{\beta+1}+O(\rho^{\beta+2}),\\
		\dot{\theta}=b\text{Cs}^{\beta+1}(\theta)\text{Sn}(\theta)\rho^{\beta}+O(\rho^{\beta+1}).
	\end{array}
\end{equation*}

We have that $\theta=0$ is a simple zero of $b\text{Cs}^{\beta+1}(\theta)\text{Sn}(\theta)\rho^{\beta}$, and thus there are solution curves of $X$ that approaching or leaving the origin with concrete slope $\theta=0$, i.e. we have characteristic orbits. The origin cannot be monodromic.

Case $\beta>n-1$ or $j^r\Phi(0)\equiv 0$: Under the generalized polar coordinates, system \eqref{eqDemoTeoAndreev2} is written as:
\begin{equation*}
	\dfrac{d\rho}{d\theta}=\dfrac{O(\rho^2)}{1+O(\rho)}.
\end{equation*}
Since the denominator does not vanish at $\rho=0$, there are no arrival directions and we can define a first return function. Thus the origin is monodromic.

Case $\beta=n-1$: System \eqref{eqDemoTeoAndreev2} can be put in the form 
\begin{equation*}
	\dfrac{d\rho}{d\theta}=\dfrac{b\rho\text{Cs}^{n-1}(\theta)\text{Sn}^2(\theta)+O(\rho^2)}{1+b\text{Cs}^{n}(\theta)\text{Sn}(\theta)+O(\rho)}.
\end{equation*}

By statement (b) of Proposition \ref{PropoGenTrigoProp}, $1+b\text{Cs}^{n}(\theta)\text{Sn}(\theta)=(\text{Cs}^{n}(\theta))^2+b\,\text{Sn}(\theta)\text{Cs}^{n}(\theta)+n\text{Sn}^2(\theta)$. Hence the discriminant of the equation $1+b\text{Cs}^{n}(\theta)\text{Sn}(\theta)=0$ is given by $(b^2-4n)\text{Sn}^2\theta$. Then, for $b^2-4n<0$, $1+b\text{Cs}^{\beta+1}(\theta)\text{Sn}(\theta)$ does not vanish, and therefore the origin is monodromic. For $b^2-4n>0$, the expression $1+b\text{Cs}^{n}(\theta)\text{Sn}(\theta)$ has simple zeros and so we have characteristic orbits. 

For $b^2-4n=0$, we make the directional blow-up $\bar{x}=x$, $\bar{y}=y/x^n$ in system \eqref{eqDemoTeoAndreev2} and conclude that the origin is a saddle-node. Therefore system \eqref{eqDemoTeoAndreev2} has characteristic orbits associated with the origin, i.e. the origin is not monodromic. This concludes the proof.\qed

\addcontentsline{toc}{chapter}{Bibliografia}
\bibliographystyle{siam}
\bibliography{Referencias.bib}
\end{document}